\documentclass{tran-l}
\usepackage{epsfig}
\usepackage[centertags]{amsmath}
\usepackage{amssymb}
\usepackage{amsthm}
\usepackage[active]{srcltx} 

\newtheorem{theo}{Theorem}[subsection]
\newtheorem{theo*}{Theorem}[subsection]
\newtheorem{coro}[theo]{Corollary}
\newtheorem{coro*}[theo*]{Corollary}
\newtheorem{lem}[theo]{Lemma}
\newtheorem{prop}[theo]{Proposition}
\theoremstyle{definition}

\theoremstyle{remark}
\newtheorem{rem}[theo]{Remark}
\numberwithin{equation}{subsection}
 
\newcommand{\eps}{\varepsilon}
\newcommand{\To}{\longrightarrow}
\newcommand{\A}{\mathcal{A}}

\newcommand{\II}{\mathcal{I}}
\newcommand{\PP}{\mathbb{P}}

\newcommand{\M}{\mathcal{M}}
\newcommand{\LL}{\cal{L}}

\newcommand{\R}{\mathbb{R}}
\newcommand{\C}{\mathbb{C}}
\newcommand{\Z}{\mathbb{Z}}
\newcommand{\T}{\mathbb{T}}
\newcommand{\Lim}{\mathop{\longrightarrow}\limits}
\newcommand{\petito}{\mathop{o(1)}\limits}
\newcommand{\cal}{\mathcal}
\newcommand{\N}{\mathbb{N}}

\newcommand{\HH}{\mathcal{H}}
\newcommand{\supp}{\mbox{supp}\,}
\def\1{{{\mbox{${\mathrm{1\negthinspace\negthinspace I}}$}}}}
\def\sqin{\sqsubset}

\begin{document}

\title[]{Gibbs measures and semi-classical approximations to action minimizing measures}

\author{ Nalini Anantharaman }

\address{Unit\'e de Math\'ematiques Pures et Appliqu\'ees,
Ecole Normale Sup\'erieure,
6, all\'ee d'Italie,
69364 LYON Cedex 07, FRANCE}

\email{nanantha@umpa.ens-lyon.fr}

\date{April 2002}


\maketitle
\section{Introduction}
\subsection{Statement of results}
Let $\R^d$ be equipped with its usual euclidean structure, and let us consider the Lagrangian
\begin{equation}\cal{L}(x,v)=\frac{\mid v\mid^2}{2}-V(x)-\langle\omega,v\rangle
\end{equation}
defined on $\R^d\times\R^d$, $V$ being a $\Z^d$-periodic function of class $C^3$, and $\omega$ an
element of $\R^d$.

This allows to define the {\em action}
$$\A(\gamma_{|[0,t]})=\int_0^t \LL(\gamma_s,\dot\gamma_s)ds$$
for any sufficiently regular path $\gamma :[0,t]\To \R^d$; for instance, $\gamma$ piecewise $C^1$,
or $\gamma$ in the space $ H_{[0,t]}=\{\gamma :[0,t]\To \R^d, \dot\gamma \in L^2([0,t])\}$.

\vspace{.2cm}
We are interested in the relations between the deterministic and stochastic dynamics, stationary in time, defined by such a system. The stationary stochastic dynamics will be described by Gibbs measures on the set of continuous paths
in $\T^d$. The deterministic dynamics will be studied through the probability
measures on $\T^d\times\R^d$, invariant under the action of the Euler-Lagrange
flow $\phi=(\phi_t)_{t\in\R}$ associated to the lagrangian (1.1.1). More precisely, the
deterministic objects corresponding to our Gibbs measures will turn out to be the action-minimizing measures in the sense of Mather; ${\cal M}_\omega\subset \T^d$ will represent the corresponding Mather set, defined as the union on $\T^d$ of the
supports of all action-minimizing measures -- more precise definitions will be given in paragraph 1.2.

The main results of the paper may be summarized as follows :

\begin{theo} Let $\HH_\beta^\omega=e^{-\beta\langle \omega,x\rangle}\circ\left(\frac{\Delta}{2\beta}
+\beta V(x)\right)\circ e^{\beta\langle \omega,x\rangle}$
and $\HH_\beta^{\omega *}=\HH_\beta^{-\omega}$
act on $C^{\infty}(\T^d)$; let $\psi_\beta,\psi_\beta^*$ be the positive eigenfunctions, associated
to their common largest eigenvalue.

Then, as $\beta\rightarrow +\infty$, the measure
$$\mu_\beta^0=\frac{\psi_\beta(x)\psi_\beta^*(x)dx}{\int_{\T^d}\psi_\beta(y)\psi_\beta^*(y)dy}$$
on $\T^d$ concentrates on the Mather set ${\cal M}_\omega$. 

If $\mu_\infty^0$ is a limit point, and $\mu_\infty$ is the (uniquely defined)
corresponding action-minimizing measure on $\T^d\times\R^d$, then $\mu_\infty$ maximizes
$$h_\phi(\mu)-\frac{1}{2}\int_{\T^d\times\R^d} (\sum_{i=1}^d\lambda_i^+(\gamma,\dot\gamma))d\mu(\gamma,\dot\gamma)$$
amongst all action minimizing measures.
\end{theo}

Here $h_\phi(\mu)$ represents the metric entropy of the invariant probability measure
$\mu$ on $\T^d\times\R^d$, with respect to the action of the Euler-Lagrange
flow $\phi=(\phi_t)_{t\in\R}$; and the
$\lambda_i^+(\gamma,\dot\gamma)$ are the $d$ first (nonnegative) Lyapunov
exponents of $(\gamma,\dot\gamma)$, under the action of $\phi$.

The theorem is valid under suitable assumptions on the Lagrangian, which will
be stated later.

Note that the measure $\psi_\beta(x)\psi_\beta^*(x)dx$ may also be written
in the form $e^{-\beta(u_\beta+v_\beta)}dx$, where $u_\beta$ is the solution
of the Hamilton-Jacobi equation with viscosity
$$-\frac{\Delta u}{2\beta}+H(x,d_x u)=C$$
($H(x,p)=\frac{\mid p+\omega\mid^2}{2}+V(x)$), and $v_\beta$ is the solution
of the same equation for the time-reversed system :
$$-\frac{\Delta v}{2\beta}+H(x,-d_x v)=C$$

\begin{coro} Let $\HH_{\hbar}=-\hbar^2\frac{\Delta}{2}+V$, and let $\psi_\hbar$
be the unique $\Z^d$-periodic positive eigenfunction, corresponding to the smallest eigenvalue of $\HH_\hbar$.

Assume that the minima of $V$ are non-degenerate.

Then, as $\hbar\To 0$, the measure
$$\frac{\psi_\hbar^2(x)dx}{\int_{\T^d}\psi_\hbar^2(y)dy}$$
concentrates to the minima of $V$ which have the smallest sum of nonnegative
Lyapunov exponents, as equilibrium points of the
flow on $\T^d\times\R^d$ associated to the differential equation
$$\ddot\gamma=+V^{\prime}(\gamma)$$
\end{coro}

We also prove the following result, as an intermediate step towards
Theorem 1.1.1 :

\begin{theo} Let $\gamma : [0,t]\longrightarrow \T^d$ be a critical point of the action
$$\A(\xi_{|[0,t]})=\int_{0}^{t}\LL(\xi_t, \dot\xi_t)dt$$ on the affine Hilbert space
$$H_{[0,t]}^{x,y}=\{\xi\in H_{[0,t]},\xi_0=x,\xi_t=y\},$$whose tangent space $H_{[0,t]}(0,0)$ is
endowed with the scalar product 
$$\langle \xi,\eta\rangle=\int_0^t \dot\xi_t.\dot\eta_t dt$$

Then the hessian $\A^{\prime\prime}(\xi)$, an autoadjoint operator on 
$H_{[0,t]}^{0,0}$, has a well defined
determinant -- the infinite product of its eigenvalues. And this determinant
coincides with the determinant of the linear endomorphism of $\R^d$, which maps $Y^\prime\in
\R^d$ to $\frac{Y_t}{t}$, where $Y_s\in T_{\gamma_s}\T^{d}$ ($s\in[0,t]$) is the solution of the linearized equation :
\begin{eqnarray*}Y^{\prime\prime}_s+V^{\prime\prime}(\gamma_s).Y_s=0\\
Y_0=0 \mbox{, } Y^{\prime}_0
=Y^{\prime}
\end{eqnarray*}
\end{theo}

\vspace{.2cm}
\subsection{More details}
Let us explain our approach to the problem.

{\bf Deterministic dynamics.} In classical mechanics, the ``principle of least action" says that
the trajectories of the system are the paths $\gamma$ which are critical points of $\A(\gamma_{|[0,t]})$, with respect
to infinitesimal variations leaving $\gamma_0$ and $\gamma_t$ fixed.
Equivalently, $\gamma$ satisfies the differential
equation \begin{equation}\ddot\gamma_s=-V^\prime(\gamma_s)\end{equation} This defines a flow on $\R^d\times\R^d$ :
$$\phi^t(x,v)=(\gamma_t,\dot\gamma_t)$$
where $\gamma$ is the solution of (1.0.2) with initial conditions $(\gamma_0,\dot\gamma_0)=(x,v)$.

If we denote $\T^d$ the $d$-torus $\R^d/\Z^d$, the flow goes to the quotient $\T^d\times\R^d$, and is called the Euler-Lagrange flow. The natural
objects for the study of the flow in the context of ergodic theory
are the invariant probability measures; that is to say, probability measures $\mu$
on $\T^d\times\R^d$ such that $\phi^t*\mu=\mu$, for all $t$. Among such measures are the {\em
action minimizing measures}, they are defined as the invariant probability measures achieving
$$\inf\{\int_{\T^d\times\R^d}\LL d\mu,\; 
\mu\mbox{ a }\phi \mbox{-invariant probability measure on }
\T^d\times\R^d\}=: -c(\omega)$$

Action minimizing measures do exist, and are supported on a compact subset of
$\T^d\times\R^d$ (\cite{Mat}).
The Mather set is then defined as
$$\tilde\cal{M}_\omega=\overline{\cup_{\mu\mbox{ \small action-min. }}\supp\mu}\subset\T^d\times\R^d$$
Action-minimizing measures are characterized by their supports : $\mu$ is action-minimizing if and only if
it is invariant and $\supp\mu\subset\tilde\cal{M}_\omega$.

The Mather set has a strong topological
property, given by Mather's graph theorem :
\begin{theo}(\cite{Mat}) The projection $\pi :\T^d\times\R^d\To \T^d$, restricted to $\tilde\cal{M}_\omega$,
is injective. Its inverse, defined on $\cal{M}_\omega=\pi(\tilde\cal{M}_\omega)$, is lipschitz.
\end{theo}
For instance, on the 2-torus ($d=2$), is $\omega$ is `close' to $0$, the Mather set will be the collection
of maxima of $V$; and for other $\omega$'s, it will be a partial lipschitz foliation of the torus.

\vspace{.3cm}
We will also be interested in the discretized model defined by
\begin{equation}L(x_0,x_1)=\frac{\mid x_1-x_0\mid^2}{2}-V(x_0)-\langle\omega,x_1-x_0\rangle\end{equation}
or, more
generally, we could consider a function $L$ with the following properties :

{\bf (Periodicity)} $L(x+n,y+n)=L(x,y)$, for all $n\in\Z^d$.

\vspace{.2cm}
{\bf (`Twist property')} For all $x\in\R^d$, $y\mapsto \partial_1 L(x,y)$
is a diffeomorphism of $\R^d$.

\vspace{.2cm}
{\bf (Superlinear growth)} $\frac{\mid L(x,y)\mid}{\mid x-y\mid}\Lim_{\mid x-y\mid\To
+\infty} +\infty$

For the discrete time system, a path $\gamma$ will simply be a sequence
$(\gamma_0,...,\gamma_n), \,\gamma_i\in\R^d$, and its action :
$$A(\gamma_0,...,\gamma_n)=\sum_{k=0}^{n-1}L(\gamma_k,\gamma_{k+1})$$
The critical points of the action (with fixed endpoints) are the paths that satisfy the recurrence relation 
$$(\gamma_{i+1}-\gamma_i)-(\gamma_i-\gamma_{i-1})=-V^\prime(\gamma_i),$$
or in the general form :
$$\partial_2 L(\gamma_{i-1},\gamma_{i})+\partial_1 L(\gamma_i,\gamma_{i+1})=0$$

This defines a diffeomorphism $\phi$ of $\R^d\times\R^d$ to itself 
$$\phi :(\gamma_0,\gamma_1)\longmapsto
(\gamma_1,\gamma_2),$$ called ``twist diffeomorphism'', which goes to the quotient $\T^d\times\R^d$.

Action minimizing measures are defined the same way, and the same results hold.

\vspace{.3cm}
{\bf Stochastic dynamics.} In the context of stochastic dynamics, the system is described
by a probability measure on the configuration space 
(the set of infinite paths) : $$W=(\R^d)^\R/\Z^d,$$ or $W=(\R^d)^\Z/\Z^d$ for the discretized
system. In this definition, the action of $\Z^d$ on $(\R^d)^\R$ or $(\R^d)^\Z$ is
given by $$(n.\gamma)(t)=\gamma(t)+n$$ for all $t\in\R$,
$n\in\Z^d$. If $[t_1,t_2]$ is an interval of time, we define similarly the
set of paths,
$$W_{[t_1,t_2]}=(\R^d)^{[t_1,t_2]}/\Z^d$$ In fact, in the
case of continuous time, all our measures will be supported on the set
of continuous paths on the torus, so that one may prefer to choose as configuration space the
set of continuous paths
$$W=C^0(\R,\R^d)/\Z^d$$ 

Once again we
will be interested in the stationary dynamics, in other words, the probability measures invariant
by translations in time :
$$\sigma^t(\gamma)(s)=\gamma(s+t)$$
Let us precise that the Borel $\sigma$-field $\cal F$ on $W$ is the smallest for which all the maps
$\gamma\mapsto \gamma_t$ are measurable; the $\sigma$-field ${\cal F}_{[t_1,t_2]}$, on $W_{[t_1,t_2]}$,
is defined the same way.

Measures on $W$ will be freely identified with measures on $C^0(\R,\R^d)$ or $(\R^d)^\Z$,
invariant under the action of $\Z^d$.

We do not require our probability measures to be carried by trajectories of the Euler-Lagrange flow;
we consider that physically significant measures are those with the {\em Gibbs property}.

For the discrete time model, a probability measure $\mu$ on $W$ will be called a Gibbs measure for
the potential $L$ if, for
all $n>0$, the conditional probability $d\mu(.|(\gamma_i)_{i\geq n},(\gamma_i)_{i\leq -n})$ is given by the expression
\begin{multline}d\mu(B|(\gamma_i)_{i\geq n},
(\gamma_i)_{i\leq -n})=\frac{\int_{(\gamma_i)_{i\in\Z}\in B} e^{-\beta\sum_{i=-n}^{n-1}L(\gamma_i,\gamma_{i+1})}
d\gamma_{-n+1}...d\gamma_{n-1}}{Z_n^\beta((\gamma_i)_{i\geq n},
(\gamma_i)_{i\leq -n})}\end{multline}
for all $B\in{\cal F}$, for $(\gamma_i)_{i\geq n}\in W_{[n, +\infty)}$ and
$(\gamma_i)_{i\leq -n}\in W_{(-\infty,-n]}$. We have denoted $d\gamma_i$ the Lebesgue measure,
and $Z_n^\beta((\gamma_i)_{i\geq n},(\gamma_i)_{i\leq -n})$ is the normalization factor which makes $d\mu(.|(\gamma_i)_{i\geq N},
(\gamma_i)_{i\leq -n})$ a probability measure.

In the expression (1.2.3), the boundary conditions $(\gamma_i)_{i\geq n}$ and
$(\gamma_i)_{i\leq -n}$ are defined modulo the action of $\Z^d$ on the set
of paths. So, given $(\gamma_i)_{i\leq -n}\in (\R^d)^{(-\infty, n]}$, the integral in (1.2.3) needs to take into account all the representatives of $(\gamma_i)_{i\geq n}$ in $(R^d)^{[n,+\infty)}/\Z^d$. As a consequence, the Gibbs measure
depends on the cohomology class $\omega$.

This definition also depends on a positive parameter $\beta$, which can be thought of,
for instance, as the inverse of a temperature.
We will denote $\mu_\beta$ the corresponding Gibbs measure : one can prove that it exists and is unique. It is ergodic
with respect to the shift $\sigma$ acting on $W$.

Note in passing that the Gibbs measure remains unchanged if
$L(x_0,x_1)$ is replaced by $L(x_0,x_1)+u(x_1)-u(x_0)+c$, where $u$ is a $\Z^d$-periodic function.

The Gibbs measure can be described as follows : there exists unique (up to a multiplicative factor) continuous positive $\Z^d$-periodic functions
$\psi_\beta$, $\psi_\beta^*$, such that
\begin{equation}\psi_\beta(x)=e^{\lambda_\beta} \int_{\R^d}e^{-\beta L(x,y)}\psi_{\beta}(y)dy\end{equation} and
$$\psi_\beta^*(x)=e^{\lambda_\beta} \int_{\R^d}e^{-\beta L(y,x)}\psi_{\beta}^*(y)dy$$
for some $\lambda_\beta\in\R$, for all $x$. Then $\mu_\beta$ is the Markov measure with initial distribution
$$\frac{\psi_\beta(x)\psi_\beta^*(x)dx}{\int_{\T^d}\psi_\beta(y)\psi_\beta^*(y)dy}$$ and with transition densities
$$P(x,dy)=\frac{\psi_{\beta}(y)}{e^{\lambda_\beta} \psi_{\beta}(x)}e^{-\beta L(x,y)}dy$$

Finally, let us mention a variational principle satisfied by $\mu_\beta$ :
amongst all $\sigma$-invariant probability measures on $W$, $\mu_\beta$ minimizes
$$\int Ld\mu-\frac{1}{\beta}H(\mu)$$
where $H$ is the functional with values in $[-\infty,0]$ defined by
$$H(\mu)=\int \log\left(\frac{d\mu(\gamma_0|\gamma_1,\gamma_2,...)}{d\gamma_0}\right)
d\mu(\gamma_0,\gamma_1,\gamma_2...)$$
if the conditional probability $d\mu(\gamma_0|\gamma_1,\gamma_2,...)$ is absolutely continuous
with respect to the Lebesgue measure; and $H(\mu)=-\infty$ otherwise.
\vspace{.2cm}

For the continuous time model, the definition of Gibbs measures is similar.
We will say that the probability $\mu_\beta$ on $W$ is a Gibbs measure for
the potential $\LL$ (given by (1.1.1)) if, for
all $T>0$, the conditional probability $d\mu(.|(\gamma_t)_{t\geq T},
(\gamma_t)_{t\leq -T})$ is given by the expression
\begin{multline}d\mu(B|(\gamma_t)_{t\geq T},
(\gamma_t)_{t\leq -T})\\=\frac{\int_{(\gamma_t)\in B} e^{\beta\int_{-T}^T V(\gamma_s)ds+\beta\langle \omega,
\gamma_T-\gamma_{-T}\rangle}
d{\cal W}_{[-T,T]}^{\beta,(\gamma_{-T},\gamma_{T})}(\gamma_{|[-T,T]})}{Z_T^\beta((\gamma_t)_{t\geq T},
(\gamma_t)_{t\leq -T})}\end{multline}
for all $B\in{\cal F}$, for $(\gamma_t)_{t\geq T}\in W_{[T, +\infty)}$ and
$(\gamma_t)_{t\leq -T}\in W_{(-\infty,-T]}$. We have denoted
$$d{\cal W}_{[-T,T]}^{\beta,(\gamma_{-T},\gamma_T)}(\gamma_{|[-T,T]})=\mbox{``}e^{-\beta\int_{-T}^T\frac{\mid\dot\gamma\mid^2}{2}}d\gamma\mbox{''}$$
the brownian bridge between $\gamma_{-T}$ and $\gamma_T$ in the time interval
$[-T,T]$, with diffusion coefficient $1/\beta$

\vspace{.2cm}
For $x,y\in\R^d$, we recall that the brownian bridge ${\cal W}^{\beta,(x,y)}_{[T_1,T_2]}$ with diffusion coefficient $1/\beta$, starting at $x$ and ending at $y$, in the time interval $[
T_1,
T_2]$
is defined as the unique
positive measure on $C^0([T_1,T_2],\R^d)$ such that
$${\cal W}^{\beta, (x,y)}_{[T_1,T_2]}(B)=(S^\beta_{t_1-T_1}\1_{B_1}S^\beta_{t_2-t_1}
\1_{B_2}
\cdots S^\beta_{t_n-t_{n-1}}(\1_{B_n}s^\beta_{T_2-t_n}(.,y)))(x)$$
for all $B\in\cal{F}_{[T_1,T_2]}$ of the form $B=\{\gamma,\gamma_{t_i}\in
B_i, \forall i\}$, where $T_1<t_1<\cdots<t_n<T_2$, and the $B_i$'s are measurable subsets of
$\R^d$. We have denoted
$$s^\beta_{t}(x,y)=\frac{1}{(2\pi t/\beta)^{d/2}} e^{-\frac{\beta\mid x-y\mid^2}{2t}} $$
the transition kernel, and
$$S_t^\beta f(x)=\frac{1}{(2\pi t/\beta)^{d/2}} \int e^{-\frac{\beta\mid x-y\mid^2}{2t}}f(y)dy,$$
for $f$ in $L^2$, the corresponding transition semi-group.
The brownian bridge ${\cal W}^{\beta,(x,y)}_{[T_1,T_2]}$
is in fact supported on the set of continuous paths with endpoints $x,y$. Note that,
in the expression (1.2.5), the boundary conditions $(\gamma_t)_{t\geq T}$ and
$(\gamma_t)_{t\leq -T}$ are defined only up to the action of $\Z^d$ on the set of paths; so that, given $(\gamma_t)_{t\leq -T}\in  C^0((-\infty, -T], \R^d)$, all the representatives of $(\gamma_t)_{t\geq T}
\in C^0([T,+\infty))/\Z^d$ need to be taken into account in the integral (1.2.5).
That is the reason why the Gibbs measure depends on the cohomology class $\omega$.

Finally, we recall that the stationary Brownian motion $\cal W$ in $\R^d$ (Wiener measure with uniform initial distribution) and the Brownian bridge may be related as follows : if $B$ is a set of the same form as before,
$${\cal W}(B)=\int \1_{B_1}(\gamma_1)...\1_{B_n}(\gamma_n){\cal W}^{\beta,\gamma_1,\gamma_2}_{[t_1,t_2]}...{\cal W}^{\beta,\gamma_{n-1},\gamma_n}_{[t_{n-1},t_n]}d\gamma_1...d\gamma_n$$
where ${\cal W}^{\beta,\gamma_1,\gamma_2}_{[t_1,t_2]}$ stands for the total mass of the measure ${\cal W}^{\beta,\gamma_1,\gamma_2}_{[t_1,t_2]}$ and has the value
$\frac{1}{(2\pi t/\beta)^{d/2}} e^{-\frac{\beta\mid \gamma_2-\gamma_1\mid^2}{2t}} $
\vspace{.2cm}

According to tastes, the parameter $1/\beta$ may be thought of as a diffusion coefficient,
a viscosity coefficient, or, as we shall explain, the $\hbar$ of quantum mechanics (but with the $\sqrt{-1}$ missing) :

There is a characterization of $\mu_\beta$ in terms of the positive eigenfunctions of the
`twisted' Hamilton operator
$$\HH_\beta^\omega=e^{-\beta\langle \omega,x\rangle}\circ\left(\frac{\Delta}{2\beta}
+\beta V(x)\right)\circ e^{\beta\langle \omega,x\rangle}$$
and its adjoint
$$\HH_\beta^{\omega *}=\HH_\beta^{-\omega}$$
which is also the twisted Hamilton operator for the time-reversed system.
Both act on $C^{\infty}(\T^d)$, and have positive eigenfunctions $\psi_\beta,\psi_\beta^*$ associated
to their common largest eigenvalue $\lambda_\beta$ :
\begin{equation}\HH_\beta^\omega \psi_\beta=\lambda_\beta \psi_\beta\end{equation}
$$\HH_\beta^{\omega *} \psi_\beta^*=\lambda_\beta \psi_\beta^*$$
Then $\mu_\beta$ is the Markov process with initial distribution
$\frac{\psi_\beta(x)\psi_\beta^*(x)dx}{\int_{\T^d}\psi_\beta(y)\psi_\beta^*(y)dy}$, and with transition semi-group
$$f\mapsto P^t_{\beta,\omega}f=\frac{1}{e^{\lambda_\beta t}\psi_\beta}\exp t\HH_\beta^\omega.(\psi_\beta f)$$

\vspace{.3cm}

{\bf Limiting behaviour as $\beta\To +\infty$.}

A natural question is to find the behaviour of the Gibbs measure $\mu_\beta$
as $\beta\To +\infty$.

We will say that a sequence of probability measures $(\mu_n)_{n\geq 0}$ on
$(W,\cal F)$ converges to $\mu$, if
$$\int f d\mu_n\To \int f d\mu,$$
for all $f$ on $W$ of the form
$$\gamma\longmapsto g(\gamma_{t_1},...,\gamma_{t_l})$$
for some $t_1<...<t_l$, and $g$ a bounded continuous function on $(\R^d)^l/\Z^d$.

\begin{lem} There exists a sequence $\beta_k \To +\infty$ such that
$(\mu_{\beta_k})_{k\geq 0}$ converges.
\end{lem}

(The proof is given in Part 2).

We would like to know whether the limit is independent of the sequence $(\beta_k)$, or not. It is not too hard to see (Corollary 2.0.15) that any limit point $\mu_\infty$ is
carried on the subset of $W$ formed by trajectories of the Euler-Lagrange
flow (or the twist diffeomorphism), so that it can be naturally
identified to a probability measure on $\T^d\times
\R^d$, invariant under the flow; and this measure is in fact action minimizing.
But there can be several action-minimizing measures,
with the same support, or with different supports (included in the Mather set), and
we would like to know which of them can appear as limits of the Gibbs measures
defined above.

The following theorem gives a partial answer, by ruling out
certain action-minimizing measures as limit points :

\begin{theo} (a) Let $\mu_\infty$ be a limit point of the family $(\mu_\beta)_{\beta\rightarrow
\infty}$; then
$\mu_\infty$ is carried by trajectories of the Euler-Lagrange flow (or the
twist diffeomorphism, in the discrete time model) associated
to $\LL$; thus, it can be identified with a measure on $\T^d\times\R^d$, invariant under the flow.

This measure is an action minimizing measure.

(b) Moreover, under
the technical assumptions (A1), (A2) and (A3) below, if $\mu$ is another action-minimizing measure, then
\begin{equation}h_\phi(\mu)-\frac{1}{2}\int_{\T^d\times\R^d} (\sum_{i=1}^d\lambda_i^+(\gamma,\dot\gamma))
d\mu(\gamma,\dot\gamma)
\leq h_\phi(\mu_\infty)-\frac{1}{2}\int_{\T^d\times\R^d} (\sum_{i=1}^d\lambda_i^+(\gamma,\dot\gamma))
d\mu_\infty(\gamma,\dot\gamma),\end{equation}
where :

-- $h_\phi(\mu)$ denotes the metric entropy of $\mu$ with respect to the action of $\phi$.

-- the $\lambda_i^+$ are the $d$ first (nonnegative) Lyapunov exponents
of $(\gamma,\dot\gamma)$ for the Euler-Lagrange flow.
\end{theo}

The $2d$ Lyapunov exponents of $(\gamma,\dot\gamma)$ under $\phi$, which are the same as those of the Hamiltonian flow (or exact symplectomorphism,
in discrete time)
obtained by Legendre duality, come into pairs $(\lambda_i^+, -\lambda_i^+)_{1\leq i\leq d}$,
with $\lambda_i^+\geq 0$.

\begin{rem}Note that the entropies $h_\phi(\mu)$ and $h_\phi(\mu_\infty)$ will vanish automatically
in the following cases :

-- $\omega=0$.

-- $d=2$, continuous time.

-- $d=1$, discrete time.

In all these cases, the theorem says that the Gibbs measures will converge to
the `least hyperbolic' action-minimizing measures. Otherwise, there will be a competition between entropy and Lyapunov exponents in order to decide of the limiting measure.
\end{rem}

\begin{rem} The theorem holds for
a mechanical Lagrangian, of the form (1.1.1) or (1.2.2). For
a more general Lagrangian, the question just does not make
sense in continuous time (because of the special part played by brownian motion);
but in discrete time and for a Lagrangian of more general
form, the Gibbs measures are still well defined,
and the conclusion of Theorem (1.2.3) has to be slightly
modified. What we get is
\begin{multline*}h_\phi(\mu)-\frac{1}{2}\int_{\T^d\times\R^d} (\sum_{i=1}^d\lambda_i^+(\gamma_0,\gamma_1)
+\log\mid \partial^2_{12} L(\gamma_0,\gamma_1)\mid)d\mu(\gamma_0,\gamma_1)\\
\leq h(\mu_\infty)-\frac{1}{2}\int_{\T^d\times\R^d} (\sum_{i=1}^d\lambda_i^+(\gamma_0,
\gamma_1)+\log\mid \partial^2_{12} L(\gamma_0,\gamma_1)\mid)
d\mu_\infty(\gamma_0,\gamma_1)\end{multline*}
instead of the simpler inequality of Theorem (1.2.3). 
\end{rem}

Let us define, once and for all, our notations for path spaces, and give the assumptions (A1), (A2), (A3) under which the theorem holds :

{\bf Paths spaces.} Let us summarize our notations for the various paths spaces we use.

As the reader will notice, we shall not make a clear distinction between a path on the torus (an element of $C^0(\R,\R^d)/\Z^d$ or $(\R^d)^\Z/\Z^d$) and a lift to $\R^d$.

We denote $H_{[0,t]}$ the Hilbert manifold of paths $[0,t]\To \T^d$, with $L^2$ derivative. The scalar product
is denoted $\langle .,.\rangle$; for $x,y\in\R^d$, $H_{[0,t]}^x$ is the affine
subspaces of paths starting at $x$, and $H_{[0,t]}^{x,y}$ the space of paths which can be lifted to a path in $R^d$ with endpoints $x,y$.

We denote $W_{[0,t]}$ the Banach manifold of continuous paths $[0,t]\To \T^d$.
The topology is that of uniform convergence on compact subintervals; $W_{[0,t]}^x$ and $W_{[0,t]}^{x,y}$ are, respectively, the affine
subspaces of paths starting at $x$, and with endpoints $x,y$.

In the continuous time model, the space $W_{[0,t]}$ can be endowed with the Wiener measure starting at $x$ ${\cal W}^x_{[0,t]}$, carried
on $W_{[0,t]}^x$, or by the brownian bridge ${\cal W}^{x,y}_{[0,t]}$, carried
on $W_{[0,t]}^{x,y}$.

\vspace{.2cm}
{\bf (A1)} For all $x,y\in\R^d$, for all $t$, the action $\A$ has only non-degenerate minima
in $H_{[0,t]}^{x,y}$,
and the number of minimizers is bounded, independently of $x,y,t$.

{\em In order to simplify the notations, we will assume in the proof
that there is only one (non-degenerate) minimizer, for all $x,y,t$. We will denote $$h_t(x,y)=\inf_{h_{[0,t]}^{x,y}}
\A$$For all $t$, $h_t$ is a lipschitz function, with lipschitz constant independent
of $t$ for $t\geq 1$.}

\vspace{.2cm}
{\bf (A2)} There exists $\eps_0 >0$ such that, for all $\eps\leq\eps_0$, for all $t$, if $\gamma_0,
\gamma_t\in\R^d$ are such that
$$\mid \gamma_0 -\xi_0\mid \leq\eps$$
$$\mid \gamma_t -\xi_t\mid \leq\eps$$
for some $\xi$ in the Mather set,
then there exists a minimizer $\bar\gamma\in H_{[0,t]}^{\gamma_0,\gamma_t}$ of
$\A :H_{[0,t]}^{\gamma_0,\gamma_t}\To\R$ 
such that
$\mid \gamma_s-\xi_s\mid\leq \eps$ for all $0\leq s\leq t$.
\vspace{.2cm}

{\bf (A3)} It is possible to replace the Lagrangian $\LL(x,v)$ (respectively $L(\gamma_0,\gamma_1)$)
by a cohomologous
Lagrangian $\LL(x,v)-d_xu.v +c$ (respectively $L(\gamma_0,\gamma_1)-u(\gamma_1)+u(\gamma_0)+c$) which is
nonnegative, and vanishes on the Aubry-Mather set (see Remark 1.0.14). And this can be done in such a way that
$$\mbox{Leb}(\{(\gamma_0,\gamma_t)\in (\R^d)^2/\Z^d, h_t(\gamma_0,\gamma_t)\leq\eps\})
\leq B(t)\eps^{d/2}$$
with $\lim_{t\To\infty} \frac{\log B(t)}{t}=0$; in other terms,
$$\beta^d\int_{(\R^d)^2/\Z^d} e^{-\beta h_t(\gamma_0,\gamma_t)}d\gamma_0d\gamma_t \leq B(t)$$

\begin{rem} Assumption (A3) is on the non-degeneracy of the Aubry-Mather set as the set
of global minimizers of the action. For instance, it is satisfied for $L(x,v)=\frac{\mid v\mid^2}{2}-V(x)$
where $V$ has only non degenerate maxima. The first part of the assumption, about the existence of $u$, is justified by a recent result
by Fathi and Siconolfi, see Remark 1.2.11.
\end{rem}

\begin{rem}
As the reader who goes through Part 2 may see, these assumptions are not the optimal ones under which
the theorem holds (however, it does not seem possible to completely get rid of them).
For instance, (A1) could probably be replaced by much weaker bounds
on the number of minimizers (which still have to be non degenerate) : it seems enough to ask for the number
of minimizers of $\A :H_{[0,t]}^{x,y}\To\R$ to grow subexponentially fast in $t$. I also found
conditions which look weaker than (A3), but not very natural.
\end{rem}

\begin{rem}
As A. Fathi pointed out to me, there is no reason a priori that there should exist a minimizing
measure achieving the variational principle (1.2.7); it seems that assumptions (A2) and (A3) will ensure this. 
\end{rem}

Theorem 1.1.1 and Corollary 1.1.2 are direct consequences of Theorem 1.2.3;
in Theorem 1.1.1, we have to assume (A1), (A2) and (A3).

Theorem 1.1.3 is an ingredient towards Theorem 1.2.3;
the similar statement, for a twist diffeomorphism generated by a function of the form (1.2.2), was
already known to a number of people. It is stated and proved in part 2.

\begin{rem} More generally, one can hope that Theorems 1.2.3 and 1.1.3 should hold for a Euler-Lagrange flow
associated to a Lagrangian of the form (1.1.1), on a compact Riemannian manifold (maybe with a modification
due to curvature). But the proof would involve even more technicalities
than in the flat case.
\end{rem}

\vspace{.3cm}

We conclude this part drawing a few connections with some existing works on
Hamilton-Jacobi equations.
 
{\bf Hamilton-Jacobi equations.}

There is a natural relation between the Gibbs measures and action-minimizing measures, and the solutions of the Hamilton-Jacobi equation (with or without
viscosity).

Let $H:\T^d\times\R^d\To \R$ be the Hamiltonian associated to the Lagrangian $L$.
More explicitely, $H(x,p)=\frac{\mid p+\omega\mid^2}{2}+V(x)$.
The Hamilton-Jacobi equation with a viscous term reads

\vspace{.2cm}
{\bf (HJV)} $-\frac{\Delta u}{2\beta}+H(x,d_x u)=C$

($\frac{1}{\beta}$ playing the role of a viscosity coefficient),
and the same without the viscosity term is the usual stationary Hamilton-Jacobi
equation :

\vspace{.2cm}
{\bf (HJ)} $H(x,d_x u)=C$

\vspace{.2cm}

The reference for the study of solutions of these equations is the book of Lions, \cite{Lio}.

Equation (HJV) only has a solution for the value $C=\frac{\lambda_\beta}{\beta}$,
and this solution is unique, given by $u_\beta=-\frac{\log \psi_\beta^*}{\beta}$. Similarly, if we
considered (HJV) for the reversed Hamiltonian $H(x,p)=\frac{\mid p-\omega\mid^2}{2}+V(x)$,
the corresponding solution would be $v_\beta=-\frac{\log \psi_\beta}{\beta}$. Thus,
the measure $\mu_\beta^0$ (the marginal at $t=0$
of the Gibbs measure) is, up to renormalization, $e^{-\beta(u_\beta(x)+v_\beta(x))}dx$.

\vspace{.2cm}
\begin{rem} This also shows that $\mu_\beta^0$ coincides with the projection on
$\T^d$ of the ``stochastic Mather measures'' studied by Gomes in \cite{Gom} (however
the Gibbs measures and stochastic Mather measures themselves are not the same objects).
\end{rem}

For the equation (HJ), the natural notion of solution is that of ``viscosity
solution'' (see \cite{Lio}). Such solutions exist only for a certain value of
$C$, which, after the works of Ma\~ne, Mather, Fathi... (\cite{Mn1}, \cite{Mn2}, \cite{Mat}, \cite{Fa2}), is $C=c(\omega)$. Some of these solutions,
possibly not unique, are lipschitz. We will denote $S_-$
the set of lipschitz viscosity solutions of (HJ), and $u_-$ an element of  $S_-$.
If we consider the Hamilton-Jacobi equation associated with the time reversed system,
we obtain a second class $S_+$ of lipschitz viscosity solutions $-u_+$.

An equivalent way of finding solutions $u_-$ or $u_+$ is as fixed points, respectively,
of the Hopf-Lax semi-groups :
$$T_t^-
u(x)=\inf_{\gamma \in C^1( [-t,0], \T^d),
\gamma(0)=x}\{u(\gamma_{-t})+\int_{-t}^0 \LL(\gamma_s,\dot\gamma_s)ds +c(\omega)t\}
$$ and of $$T_t^+
u(x)=\sup_{\gamma \in C^1( [0,t], \T^d),
\gamma(0)=x}\{u(\gamma_t)-\int_0^t \LL(\gamma_s,\dot\gamma_s)ds -c(\omega)t\}
$$

The elements of $S_-, S_+$ come naturally into pairs $(u_-, u_+)$, called conjugate solutions,
satisfying $u_--u_+=0$ on the Mather set ${\cal M}_\omega$, and $u_--u_+\geq 0$ elsewhere.

The graphs $\{(x,d_x u_-)\}$, $\{(x,d_x u_-)\}$, when transported by Legendre duality
to the tangent space, are respectively invariant by $(\phi_t)_{t\leq 0}$, and $(\phi_t)_{t\geq 0}$.
The intersection of these two sets is a $(\phi_t)$-invariant subset of $\T^d\times\R^d$,
denoted $\tilde \II_{(u_-,u_+)}$;
it contains the Mather set and has the same graph property (Theorem 1.2.1), but may contain, in addition,
orbits which do not lie in the support of an action-minimizing measure. Its projection to $\T^d$, $\II_{(u_-,u_+)}$,
is the set of points where $u_--u_+=0$.
The set $\cup_{(u_-,u_+)} \tilde \II_{(u_-,u_+)}$ is called the Ma\~ne set, and
$\cap_{(u_-,u_+)} \tilde \II_{(u_-,u_+)}$  is called the Aubry set.
One can show that the Ma\~ne set is the set of ``globally'' action minimizing trajectories,
and that the Aubry set is, roughly, the accumulation points of closed curves which are `almost'
action minimizing (for more details, see the work of Fathi, \cite{Fa1}, \cite{Fa2}).

For the discretized system, the same results hold with the fixed points
$u_-, u_+$ of the (nonlinear) operators :
$$T^- u(x)=\inf_y \{u(y)+L(y,x)+c(\omega)\}$$
and
$$T^+ u(x)=\sup_y\{u(y)-L(x,y)-c(\omega)\}$$

For the Hamilton-Jacobi equation (HJV), the behaviour of the solution
$-\frac{\log\psi_\beta^*}{\beta}$ as $\beta\To +\infty$, is already a subject of great interest. The
family $(-\frac{\log\psi_\beta^*}{\beta})_{\beta >0}$ can be shown to be uniformly lipschitz, and any limit
point (in uniform topology) will
be a viscosity solution $u_-$ of (HJ) without viscosity. In particular, this implies that $$\frac{\log
\lambda_\beta}{\beta}\Lim_{\beta\To +\infty} c(\omega)$$

The problem of the existence of a (unique) limit of $(-\frac{\log\psi_\beta^*}{\beta})_{\beta\rightarrow \infty}$,
has been studied in \cite{JKM}, in the particular case $d=1$ and $\omega=0$. This question,
although obviously related to our problem, is not exactly of the same nature. 
The existence of a limit for $-\frac{\log\psi_\beta^*}{\beta}$
and $-\frac{\log\psi_\beta}{\beta}$ yields a Large Deviation property of the
family $(\mu_\beta)$, whereas we are interested in the existence of a weak limit. The large
deviation property influences the possible choice of a weak limit, and vice-versa,
but the two phenomena are not equivalent. Fortunately, the result of \cite{JKM}
is compatible with ours !

\begin{rem} We note that, if $u$ is an element of $S_-$ or $S_+$, we can replace the action $\A
(\gamma_{|[0,t]})$ by $\A
(\gamma_{|[0,t]})-u(\gamma_t)+u(\gamma_0) +c(\omega)t$.
This way, the action of a path is always nonnegative, and the action
of a trajectory in the Mather set is zero. This transformation does not
change the notion of Gibbs measure, nor the derivatives of the action, for
fixed endpoints.

Fathi and Siconolfi have a recent result according to which $\LL -du+c(\omega)$ can actually
be made non-negative, and vanishing {\em precisely} on the Aubry set, for some $u$
of class $C^1$(this smooth function $u$
will, a priori, not belong to $S_-$ nor $S_+$). This legitimates Assumption (A3).

So, if we add to $\LL$ a constant and an exact $1$-form (which does not change the Gibbs measures),
we can assume in the rest of the paper that $\LL\geq 0$, and vanishes precisely on the Aubry set.
\end{rem}

\section{Proof of the results}

{\em Proof of Lemma 1.2.2.}

We give the proof in the case of continuous time, the
case of discrete time is similar but requires less arguments.

Let us fix $T>0$. To get rid of some constants, assume 
that $\parallel\omega\parallel\leq 1$ and $\mid V\mid\leq 1$. For all $0< t\leq T$, for all $x\in\R^d$,
\begin{multline*}\mu_\beta(\mid \gamma_t-\gamma_0\mid\geq 4dt\;|
\gamma_0=x)=\frac{\int_{W_{[0,t]}} \1_{\{\mid \gamma_t-x\mid\geq 4dt\}}
e^{\beta\int_0^t V(\gamma_s)ds +\beta\langle\omega,\gamma_t-x\rangle}
d{\cal W}^{\beta,x}_{[0,t]}(\gamma)}{
\int_{W_{[0,t]}} 
e^{\beta\int_0^t V(\gamma_s)ds +\beta\langle\omega,\gamma_t-x\rangle}
d{\cal W}^{\beta,x}_{[0,t]}(\gamma)}\\
\leq \frac{\int_{W_{[0,t]}} \1_{\{\mid \gamma_t-x\mid\geq 4dt\}}
e^{\beta( t +\mid\gamma_t-x\mid)}
d{\cal W}^{\beta,x}_{[0,t]}(\gamma)}{
\int_{W_{[0,t]}} e^{-\beta( t +\mid\gamma_t-x\mid)}
d{\cal W}^{\beta,x}_{[0,t]}(\gamma)}\\
= \frac{e^{2\beta t}\int_{\R^d}\1_{\{\mid y\mid\geq 4dt\}}
e^{\beta\mid y\mid-\beta\frac{\mid y\mid^2}{2t}}dy}{\int_{\R^d}
e^{-\beta\mid y\mid-\beta\frac{\mid y\mid^2}{2t}}dy}\\
\leq \frac{e^{2\beta t}\int_{\R^d}\1_{\{\mid y\mid\geq 4t\}}
e^{-\beta\frac{\mid y\mid^2}{4dt}}dy}{\int_{\R^d}
e^{-\beta(\mid y\mid+\frac{\mid y\mid^2}{2t})}dy}\\
\lesssim \mbox{Cst }e^{-2\beta t}\beta^{d/2}
\end{multline*}
for all $t>0$ and $\beta$ large enough; we have used the following estimate
for Brownian motion im $\R^d$ :
$$\PP(\frac{\gamma_t}{\sqrt{\beta}}\geq \delta)\leq 4de^{-\frac{\beta\delta^2}{4dt}}$$
(cf \cite{DZ}, (5.2.2)).

As a consequence, for all $t\not=s\leq T$,
\begin{equation}\mu_\beta(\mid \gamma_t-\gamma_s\mid\geq 4\mid t-s\mid)\lesssim \mbox{Cst }e^{-2\beta \mid t-s\mid}
\beta^{d/2}\end{equation}

This implies in particular the tightness of the laws of $\gamma_t$ under $(\mu_\beta)_{\beta>0}$, for all $t$;
so that we can find a subsequence $\beta_k \To +\infty$ such that
$$\mu_{\beta_k}(g(\gamma_{t_1},...,\gamma_{t_l}))\Lim_{k\To +\infty}
\mu_{\infty}(g(\gamma_{t_1},...,\gamma_{t_l}))$$
for some $\mu_\infty$, if
$t_1<...<t_l$ range over a dense denumerable subset of $[0,T]$, and
$g$ is a bounded continuous function on $(\R^d)^l/\Z^d$.

But actually, thanks to inequality (2.0.8), the convergence will take place for all
$t_1<...<t_l\in [0,T]$, and
$g$ bounded continuous function on $(\R^d)^l/\Z^d$.

\begin{prop} (a) Let $\psi_\beta,\psi_\beta^*$ be as in (1.2.4) or (1.2.6). Then
the families of functions $(-\frac{1}{\beta}\log\psi_\beta)_{\beta >0},
(-\frac{1}{\beta}\log\psi_\beta^*)_{\beta >0}$ are equilipschitz.

(b) If $\beta_k\To +\infty$ is a sequence such that
$$-\frac{1}{\beta_k}\log\psi_{\beta_k}\To -u_+$$
and
$$-\frac{1}{\beta_k}\log\psi^*_{\beta_k}\To v_-$$
in the uniform topology,
for some continuous functions $u_+$ and $v_-$, then $u_+\in S_+$ and $v_-\in S_-$.

(c) Let $J=\inf (v_--u_+)$, so that $$-\frac{\log\psi_{\beta_k}+\log\psi_{\beta_k}^*}{\beta_k}
+\frac{\log\int \psi_{\beta_k}(y)\psi_{\beta_k}^*(y)dy}{\beta_k}\Lim v_- -u_+-J,$$
and let $u_-$ be the function in $S_-$ conjugate to $u_+$,
then $u_-\leq v_- -J$.
\end{prop}

\begin{proof} The first assertions are well known results about viscosity
solutions of (HJV), and the vanishing viscosity method (\cite{Ba}).

As to the last assertion, it is a consequence of the inequality $v_--u_+ -J\geq 0$,
and the characterizations of the conjugate solution $u_-$ as the smallest element in $S_-$
satisfying $u_--u_+\geq 0$ (\cite{Fa2}). 
\end{proof}

\begin{prop}(Large deviation upper bound) Let $t>0$. Then for any
subset $K\subset W_{[0,t]}$, closed for the uniform topology,
$$\limsup \frac{1}{\beta}\mu_\beta(K) \leq - \inf_{\gamma\in K} \inf_{(u_-,
u_+)}u_-(\gamma_0)+ \A(\gamma_{|[0,t]})-u_+(\gamma_t)+tc(\omega)$$
where the first $\sup$ is taken over the set of conjugate fixed points of the Hopf-Lax semi-groups.
\end{prop}

\begin{rem} As mentioned in the first part, for $d=1$ and $\omega=0$
a sufficient condition of existence of a large deviation principle (with upper and
lower bounds) is given in \cite{JKM}.
\end{rem}

\begin{coro} If $\mu_\infty$ is a limit point of $\mu_\beta$, it is
carried by trajectories of the Euler-Lagrange flow (or twist diffeomorphism), and corresponds
to an action-minimizing measure on $\T^d\times\R^d$.
\end{coro}

\begin{proof} {\em (Corollary 2.0.15)} After Proposition 2.0.13,
the measure of a closed set $K\subset W$ will go to zero exponentially fast, unless $K$ contains
trajectories $\gamma$, such that
$$\inf_{(u_-,u_+)}u_-(\gamma_0)+ \A(\gamma_{|[0,t]})-u_+(\gamma_t)+tc(\omega)
=0,$$
$t$ arbitrarily large. In other words, $K$ must intersect the Ma\~ne set.

But all the invariant measures carried by the Ma\~ne set are, in fact, carried by the
Mather set, and action-minimizing.
\end{proof}

\begin{proof}{\em (Proposition 2.0.13)}
Recall the expression of $\mu_\beta(K)$, for $K\subset W_{[0,t]}$ :
\begin{multline*}\mu_\beta(K)=\frac{e^{-t\lambda_\beta}}{\int_{\T^d}
\psi_\beta \psi_\beta^*}\int_{\gamma_0\in\T^d}
\psi_\beta^*(\gamma_0)d\gamma_0\left(\int_{\gamma\in K}e^{\beta\int_0^t V(\gamma_s)ds
+\beta\langle\omega,\gamma_t-\gamma_0\rangle}\psi_\beta(\gamma_t)
d{\cal W}^{\beta,\gamma_0}_{[0,t]}(\gamma)\right)
\end{multline*}
We have already seen that $\frac{\lambda_\beta}{\beta}\Lim_{\beta\To +\infty}c(\omega)$.
We also recall that, for all $x\in\T^d$,
\begin{multline}\limsup\frac{1}{\beta}\log 
\int_{\gamma\in K}e^{\beta\int_0^t V(\gamma_s)ds
+\beta\langle\omega,\gamma_t-\gamma_0\rangle}e^{-\beta u(\gamma_t)}
d{\cal W}^{\beta,x}_{[0,t]}(\gamma)
\leq -\inf_{\gamma\in K,\gamma_0=x}\A(\gamma_{|[0,t]})+u(\gamma_t)
\end{multline}
for every continuous function $u$ on $\T^d$, from the large deviation
results of Schilder and Varadhan (\cite{Schi}, \cite{Var}, \cite{DZ}).

Finally, let us consider a sequence $\beta_k\rightarrow +\infty$ such that $\frac{1}{\beta_k}\log\mu_{\beta_k}(K)$
converges in $\R\cup\{-\infty\}$. Keeping the notations of Proposition 2.0.12, we
may also assume that $$-\frac{1}{\beta_k}\log\psi_{\beta_k}\To -u_+\in S_+$$
$$-\frac{1}{\beta_k}\log\psi^*_{\beta_k}\To v_-\in S_-$$
and
$$\frac{1}{\beta_k}\log\int_{\T^d}
\psi_\beta \psi_\beta^*\To -J,$$
with $v_--J$ larger than the function $u_-$ conjugate to $u_+$.

Combining this with (2.0.9), we get
\begin{multline*}
\limsup \frac{1}{\beta_k}\log\mu_{\beta_k}(K) \leq - \inf_{\gamma\in K} v_-(\gamma_0)+
\A(\gamma_{|[0,t]})-u_+(\gamma_t)+tc(\omega)-J
\\ \leq - \inf_{\gamma\in K} u_-(\gamma_0)+ \A(\gamma_{|[0,t]})-u_+(\gamma_t)+tc(\omega)
\\ \leq - \inf_{\gamma\in K} \inf_{(u_-, u_+)}u_-(\gamma_0)+ \A(\gamma_{|[0,t]})-u_+(\gamma_t)+tc(\omega)
\end{multline*}
Since this is true for every subsequence $\beta_k$, we have proved Proposition
2.0.13.
\end{proof}

\subsection{Proof for the discrete time model}
We now turn to the proof of Theorem 1.0.3, in the discrete time case. We tried to choose arguments which are transposable to the case of continuous time.

In discrete time, the $(n-1)d$-dimensional path space $H_{[0,n]}^{0,0}=W_{[0,n]}^{0,0}$ can be endowed with a large choice of euclidean structures,
and we choose the simplest :
$$(\gamma,\gamma)=\sum_{i=1}^{n-1} \parallel \gamma_i\parallel^2$$ (but
note that, when passing to continuous time, we will need to use $\sum_{i=0}^{n-1}\parallel \gamma_{i+1}
-\gamma_i\parallel^2$ instead).

Let $A^{\prime\prime}(\gamma)$ be the hessian matrix at $
\gamma\in W$ of the (formal) sum $A(\gamma)=\sum_{k\in\Z}L(\gamma_k,\gamma_{k+1})$.
We see $A^{\prime\prime}(\gamma)$ as an infinite symmetric matrix, which can be decomposed into $d\times d$ blocks $(A^{\prime\prime}_{ij})_{i,j\in\Z}$ :
$$A^{\prime\prime}_{ii}= \partial^2_{11}L(\gamma_{i-1}, \gamma_i)+\partial^2_{22}L(\gamma_i, \gamma_{i+1})$$
and
$$A^{\prime\prime}_{i, i+1}=\partial_{21} L(\gamma_i, \gamma_{i+1})$$
This way, the $nd\times nd$ submatrix $_nA^{\prime\prime}(\gamma)$, corresponding
to indices $1\leq i,j\leq n$, is the hessian matrix of the action
$A(\gamma_{|[0,n+1]})$
with respect to the variables $\gamma_1,\cdots,\gamma_n$.
\vspace{.3cm}

{\bf Notation :} -- In what follows, we shall denote $[M]$ the determinant of a square
matrix $M$.

-- unless stated otherwise, we shall always represent matrices in $d$-block
form; for instance, if $M$ is an $nd\times nd$ matrix, $M_{ij}$ ($1\leq i,j\leq n$) will be the $d\times d$ block in position $(i,j)$.

-- if $\gamma_0,\gamma_n\in\R^d$, we will denote
$_nA^{\prime\prime}(\gamma_0,\gamma_n)$ the hessian of the action $A :
H^{\gamma_0,\gamma_n}_{[0,n]}\To\R$ at its minimizer (which has been assumed unique
for simplicity). If $\gamma$ is a minimizer, then $_nA^{\prime\prime}(\gamma_0,\gamma_n)=\;_nA^{\prime\prime}(\gamma)$.

\vspace{.2cm}
The following theorem will be the first step towards Theorem 1.2.3 :
\begin{theo} Let $\mu$ be an action-minimizing measure, and $\mu_\infty$ 
a limit point of $(\mu_\beta)_{\beta\rightarrow +\infty}$. Then, under the assumptions (A1), (A2) and (A3),
$$h_\phi(\mu)-\frac{1}{2}\int_W \lim_n \frac{1}{n}\log[_nA^{\prime\prime}(\gamma)]d\mu(\gamma)
\leq h_\phi(\mu_\infty)-\frac{1}{2}\int_W
\lim_n \frac{1}{n}\log[_nA^{\prime\prime}(\gamma)]d\mu_\infty(\gamma)$$
\end{theo}

(In the notations we will not distinguish $\mu$ and $\mu_\infty$, $\sigma$-invariant probability measures on $W$
carried by trajectories
of $\phi$,
from the $\phi$-invariant probability measures on $\T^d\times\R^d$ which naturally correspond to them.)

The second step will be the following relation between determinants and Lyapunov exponents (found in
a paper by Thouless, \cite{Thou}) :
\begin{prop} If $\mu$ is an action-minimizing measure on $W$, then the limit $\lim
\frac{1}{n}\log[_nA^{\prime\prime}(\gamma)]$ exists for $\mu$-almost every $\gamma$, and is equal to
$$\sum_1^d \lambda_i^+(\gamma),$$
the sum of the $d$-first (nonnegative) Lyapunov exponents of $(\gamma_0,\gamma_1)$
under the twist diffeomorphism $\phi$.
\end{prop}

In order to prove Proposition 2.1.2, we will need the following facts, obtained
by basic manipulations of determinants of symmetric matrices :
\begin{lem} Let $M$ be a symmetric matrix, decomposed in the form 
\begin{displaymath}
M=\left(
\begin{array}{ccc}A & ^tC\\
C &B
\end{array}\right)
\end{displaymath}
(where $A$ and $B$
are square symmetric matrices, and $C$ is a rectangular matrix of appropriate dimension).

Then $[M]=[A].[B-CA^{-1}\;^tC]$.

If $M$ is (definite) positive, then $A$ and $B-CA^{-1}\;^tC$ are (definite) positive, and
$$[M]\leq[A].[B]$$
\end{lem}

\begin{proof} {\em (Proposition 2.1.2)}

Lemma 2.1.3 implies a property of subadditivity of $\log[_nA^{\prime\prime}(\gamma)]$ :
\begin{lem} If $\gamma\in W$ is such that $(\gamma_0,\gamma_1,...,\gamma_{n+1})$ is a minimizer
of the action with fixed endpoints, then, for all $m\leq n$,
$$[A^{\prime\prime}(\gamma)]\leq [A^{\prime\prime}_{m}(\gamma)].[A^{\prime\prime}_{n-m}(\sigma^m
\gamma)]$$
\end{lem}
According to the subbaditive ergodic theorem (\cite{Kin}), this implies the existence of $\lim
\frac{1}{n}\log[_nA^{\prime\prime}(\gamma)]$ for $\mu$-almost every $\gamma$, if $\mu$ is action-minimizing. Let us
now identify this limit with the Lyapunov exponents.

\begin{lem} Let $(\gamma_i)_{0\leq i\leq n}$ be a trajectory of the twist diffeomorphism
$\phi$. Let us consider the equation of variations, along $(\gamma_i)$ :
$$(Y_{i+1}-Y_i)-(Y_i-Y_{i-1})+V^{\prime\prime}_{\gamma_i}.Y_i=0$$
with an initial condition $Y_0=0$.

Then, for all $n$, the determinant of the linear map $Y_1\longmapsto Y_n$ (from $\R^d$ to $\R^d$) is equal to the determinant
of the $(n-1)d\times (n-1)d$ matrix $_{n-1}A^{\prime\prime}(\gamma)$.
\end{lem}

\begin{proof} {\em (Lemma 2.1.5).}

Let us assume that $_n A^{\prime\prime}(\gamma)$ is invertible. Me may then decompose the matrix $G=\;_nG={_nA^{\prime\prime}(\gamma)}^{-1}$
into $d\times d$ blocks $(G_{ij})_{1\leq i,j\leq n}$. 
A vector $Y=(Y_1,\cdots,Y_n)$ ($Y_i\in \R^d$)
satisfies $_nA^{\prime\prime}.Y=(0,0,0,\cdots,0,*)$, if and only if $Y$
is the solution of the linearized equation (?) with $Y_0=0$.
 
Equivalently, $$(Y_{n-1}, Y_{n})=d(\phi^{n-1})_{(\gamma_0,
\gamma_1)}.(0,Y_1)$$
Besides, the components
$Y_1$ and $Y_n$ are related by :
$$Y_n=G_{nn}.G_{n1}^{-1}Y_1$$

(If $_n A^{\prime\prime}(\gamma)$ were not invertible, we
could replace this expression by the well defined expression
$$Y_n=\mbox{com}(\,_nA^{\prime\prime}_{nn}).\mbox{com}(\,_nA^{\prime\prime})^{-1}_{n1}Y_1,$$ where com denotes the comatrix).

Let us evaluate the determinant of $G_{nn}.G_{n1}^{-1}$ in terms of the determinant of $A^{\prime\prime}(\gamma)$.
We first define a sequence of $d\times d$ matrices $(a_0, a_1, \cdots,a_{n-1})$ by $a_0=Id$ and
$$a_{k}=-A^{\prime\prime}_{k+1, k}(A^{\prime\prime}_{kk}+a_{k-1}A^{\prime\prime}_{k-1, k})^{-1},$$
agreeing temporarily that $A^{\prime\prime}_{01}=0$ (the sequence is well defined if
$_nA^{\prime\prime}$ has been assumed invertible).

We also define an $nd\times nd$ matrix $T$
decomposed into $d\times d$ blocks
$(T_{ij})_{1\leq i,j\leq n}$ with
$$T_{ii}=Id$$
$$T_{ij}=\prod_{k=j}^i a_{i-k}$$
(this way, $T$ is lower block-triangular). In fact, the matrix $T$ is constructed in such a way that
$D=T\;_nA^{\prime\prime}$ is an upper block triangular matrix, with blocks on the diagonal
$$D_{kk}=D_k=A^{\prime\prime}_{kk}+a_{k-1}A^{\prime\prime}_{k-1, k}$$
We have $G=D^{-1}T$ which yields immediately
$G_{nn}G_{n1}^{-1}=D_n T_{n1}^{-1}D_n^{-1}$
so that 
\begin{eqnarray*}[G_{nn}G_{n1}^{-1}]&=&[T_{n1}]^{-1}\\&=&(\prod_{k=1}^{n-1}[a_{n-k}])^{-1}
\\&=&(-1)^{nd}(\prod_{k=1}^{n}[A^{\prime\prime}_{k+1, k}])^{-1}\times\prod_{k=1}^{n-1}[D_{k}]
\\&=&(-1)^{nd}(\prod_{k=1}^{n}[A^{\prime\prime}_{k+1, k}])^{-1}\times [\;_{n-1}A^{\prime\prime}]
\end{eqnarray*}
where the last equality comes from the observation that $[\;_{n-1}A^{\prime\prime}]=[\;_{n-1}D]$. This expression is still valid even when $_nA^{\prime\prime}$ is not invertible.

Thus, the determinant of $Y_1\mapsto Y_n$ is equal to $(-1)^{nd}(\prod_{k=1}^{n}[
A^{\prime\prime}_{k+1, k}])^{-1}\times [\;_{n-1}\!A^{\prime\prime}]$.

Applying the Birkhoff and Oseledets theorems, this implies that
$$\lim \frac{1}{n}\log [\;_nA^{\prime\prime}(\gamma)]=\lambda_{(0,\R^d)}(\gamma)+\lim \frac{1}{n}\sum_{i=0}^{n-1}
\log\mid \partial^2_{12} L(\gamma_i, \gamma_{i+1})\mid$$
for $\mu$-almost every $\gamma$. Here $\lambda_{(0,\R^d)}(\gamma)$ is the Lyapunov
exponent of the subspace $(0,\R^d)$ of the tangent space $T_{(\gamma_0, \gamma_1)}(\R^d\times\R^d)$, for the
action of the diffeomorphism $\phi$
acting on $\Lambda^d(\R^d\times\R^d)$.

On the other hand,
the a.e limit
$\lim \frac{1}{n}\log [_nA^{\prime\prime}(\gamma)]$
is $\sigma$-invariant; and so must be $\lambda_{(0,\R^d)}(\gamma)$. 
Since the subspace $(0,\R^d)$ tangent at $(\gamma_0,\gamma_1)$ and the subspace 
$(0,\R^d)$ tangent at $(\gamma_1,\gamma_2)$ generate under the action of $d\phi$ the whole space of
tangent trajectories along $\gamma$, we necessarily have $$\lambda_{(0,\R^d)}(\gamma)=\sum_{i=1}^d\lambda_i^+(\gamma),$$ almost
everywhere. 
\end{proof}

We now turn to the proof of Theorem 2.1.1.
\vspace{.3cm}

{\bf A few more notations :}
-- We recall that $h_n(\gamma_0,\gamma_n)$ denotes the value of the minimum
of the action on $H^{\gamma_0,\gamma_n}_{[0,n]}$. If $\gamma\in W$, we will
denote $h_n(\gamma)=h_n(\gamma_0,\gamma_n)$.

-- if $B\subset W$, we will denote, quite informally, $(\gamma_0,\gamma_n)
\sqin B$ to say that there exists $\xi\in B$ such that $\gamma_0=\xi_0,
\gamma_n=\xi_n$.

-- if $B\subset W$, we will denote $B^\eps$ the uniform $\eps$-neighbourhood of
$B$ : $\{\gamma,\exists \xi\in B, \mid\gamma_k-\xi_k\mid <\eps, \mbox{ for
all }k\}$.
\vspace{.3cm}

\begin{proof} For simplicity we take $d=1$ in this proof.

Let $\eps >0$ and $M>0$, and consider the following (denumerable) partition of $\R^2=\R^{2d}$ :
$$\R^2=\sqcup_{i,j}\tilde P_{ij}$$
where the union runs over $i,j\in\Z, \mid j-i\mid <\frac{M}{\eps} \mbox{, or }i\in\Z, j=\infty$,
and
$$\tilde P_{ij}=\{(\gamma_0,\gamma_1), \gamma_0\in [i\eps, (i+1)\eps), \gamma_1\in [j\eps, (j+1)\eps)\}$$ 
for $\mid j-i\mid <\frac{M}{\eps}$, and
$$\tilde P_{i\infty}=\{(\gamma_0,\gamma_1), \gamma_0\in [i\eps, (i+1)\eps), \exists j,\mid j-i\mid\geq \frac{M}{\eps},
\gamma_1\in [j\eps, (j+1)\eps)\}$$
If $\eps$ is the inverse of an integer, this gives a finite partition of the quotient $W_{[0,1]}=\R^2/\Z$, and hence a finite
partition of 
$W=\sqcup P_{ij}$ :
$$P_{ij}=\{\gamma\in W, (\gamma_0,\gamma_1)\in\tilde P_{ij}\}$$
The number $M$ will be fixed later -- sufficiently large, whereas $\eps$ is doomed to tend to $0$. 

The choice of the partition $P$ induces a symbolic dynamics over a subshift in the finite alphabet $\{ P_{ij}\}$ :
$$W^P=\{(\alpha_k)_{k\in\Z}\subset\{(ij)\}^{\Z}, P_{\alpha_k}\cap\sigma^{-1}P_{\alpha_{k+1}}\not= \emptyset\}$$ If $\mu$ is a $\sigma$-invariant
measure on $W$ we will denote
$\mu^P$ its image on $W^P$. 

\vspace{.2cm}
Recall the following convexity inequality :
\begin{equation}-\sum p_i\log p_i +\sum p_i\log q_i\leq 0\end{equation}
whenever $(p_i)$ and $(q_i)$ are probability weights.

Hence, for all $n$,
\begin{multline*}-\sum_\alpha \mu(P_{\alpha_0}\cap...\cap \sigma^{-n+1} P_{\alpha_{n-1}})\log \mu(P_{\alpha_0}\cap...\cap
\sigma^{-n+1} P_{\alpha_{n-1}})\\+\sum_\alpha \mu(P_{\alpha_0}\cap...\cap \sigma^{-n+1} P_{\alpha_{n-1}})
\log \mu_\beta(P_{\alpha_0}\cap...\cap \sigma^{-n+1} P_{\alpha_{n-1}})\leq 0
\end{multline*}
the sums running over all word of length $n$ in $W^P$.

\vspace{.2cm}
{\em From now on, we will replace the $\cap$ by dots $.$ in expressions
of the type $P_{\alpha_0}\cap...\cap \sigma^{-n+1} P_{\alpha_{n-1}}$.}
 
We can rewrite this :
\begin{multline}
-\sum \mu(P_{\alpha_0}... \sigma^{-n-1} P_{\alpha_{n-1}})\log \mu(P_{\alpha_0}
..
\sigma^{-n+1} P_{\alpha_{n-1}})\\+\sum \mu(P_{\alpha_0}.. \sigma^{-n} P_{\alpha_n})\log \left(\frac{\beta}{2\pi}
\right)^{\frac{n+1}{2}}\int_{P_{\alpha_0}..
\sigma^{-n+1} P_{\alpha_{n-1}}}\!\!\!
\psi_\beta^*(\gamma_0)e^{-\beta \sum_{i=0}^{n-1}L(\gamma_i,\gamma_{i+1})}\psi_\beta(\gamma_{n})
d\gamma_0..d\gamma_{n}\\ \leq
-\sum \mu_\beta(P_{\alpha_0}..\sigma^{-n+1} P_{\alpha_{n-1}})\log \mu_\beta(P_{\alpha_0}..\sigma^{-n+1}
P_{\alpha_{n-1}})+\\
\sum \mu_\beta(P_{\alpha}
..\sigma^{-n+1} P_{\alpha_{n-1}})\log\left(\frac{\beta}{2\pi}
\right)^{\frac{n+1}{2}} \int_{P_{\alpha_0}..\sigma^{-n+1} P_{\alpha_{n-1}}}\!\!\!
\psi_\beta^*(\gamma_0)e^{-\beta \sum_{i=0}^{n-1}L(\gamma_i,\gamma_{i+1})}\psi_\beta(\gamma_{n})
d\gamma_0..d\gamma_{n}
\end{multline}
The plan is to divide by $n$, and first let $n$ tend to $\infty$; then, let $\beta$ tend to $\infty$, and then $\eps$ to $0$. 

We begin with finding an upper bound for the right hand side of the inequality, in terms of the determinants
$[A^{\prime\prime}]$. The integer $N$ is fixed and we take $n=kN$ in the inequality above.

\begin{lem} (a) {\bf (Laplace method).} Let $\gamma_0, \gamma_N\in \R^2$. Then, under Assumption (A1),
\begin{multline*}\left(\frac{\beta}{2\pi}
\right)^{\frac{N-1}{2}}\int_{\R^{N-1}}e^{-\beta\sum_{i=0}^{N-1}L(\gamma_i,\gamma_{i+1})}d\gamma_1...d\gamma_{N-1}
=\frac{e^{-\beta h_N(\gamma_0,\gamma_N)}}{[\;_{N-1}\! A^{\prime\prime}(\gamma_0,
\gamma_N)]^{1/2}}(1+\petito_{\beta\rightarrow\infty})\\
\leq \frac{1}{[\;_{N-1}\! A^{\prime\prime}(\gamma_0,
\gamma_N)]^{1/2}}(1+\petito_{\beta\rightarrow\infty})
\end{multline*}
with $\petito_{\beta\rightarrow\infty}$ uniform on each set
$\{\mid \gamma_N-\gamma_0\mid\leq K\}$.

(b) If the constant $M$, involved in the construction of the partition $P$,
is chosen large enough, 
then, for all $\gamma_0 \in \R$,
$$\left(\frac{\beta}{2\pi}
\right)^{\frac{N}{2}}\int_{(\gamma_0,\gamma_1,..,\gamma_N)\in P_{\alpha_0}..\sigma^{-N+1}P_{\alpha_{N-1}}}\!\!\!\!\!
e^{-\beta\sum_{i=0}^{N-1}L(\gamma_i,\gamma_{i+1})}d\gamma_1..d\gamma_{N}
\leq \left(\frac{\beta}{2\pi}
\right)^{\frac{N}{2}}e^{-\beta M}\leq 1$$
for $\beta$ large enough, as soon as one the $\alpha_k$'s is of the form
$i\infty$.

(c) If the constant $M$, involved in the construction of the partition $P$, is chosen large enough, 
then, for all $\gamma_0\in \T,\gamma_N\in\R$,
$$\left(\frac{\beta}{2\pi}
\right)^{\frac{N-1}{2}}\int_{(\gamma_0,\gamma_1,..,\gamma_N)\in P_{\alpha_0}..\sigma^{-N+1}P_{\alpha_{N-1}}}
e^{-\beta\sum_{i=0}^{N-1}L(\gamma_i,\gamma_{i+1})}d\gamma_1...d\gamma_{N-1}
\leq \left(\frac{\beta}{2\pi}
\right)^{\frac{N-1}{2}}e^{-\beta M}\leq 1$$
for $\beta$ large enough, as soon as one the $\alpha_k$'s is of the form
$i\infty$.
\end{lem}

Assertion (a) is the usual Laplace method, and requires Assumption (A1).
For (b) or (c), take $M$ such that $\mid\gamma_1-\gamma_0\mid >M\Rightarrow L(\gamma_0,\gamma_1)\geq\mid \gamma_1-\gamma_0\mid$; and also
use the fact that $L\geq 0$ everywhere.

We define the functions $F_N$ and $G_N^\beta$ on the subshift generated by $P$, depending
on $N$ coordinates :
$$F_{N}(\alpha_{0},..,\alpha_{N-1})=1$$
if one of the $\alpha_j$'s is of the form $i\infty$, and
$$F_{N}(\alpha_{0},..,\alpha_{N-1})=\sup\{\frac{1}{[_{N-1}A^{\prime\prime}(\gamma_0,\gamma_N)]^{1/2}},
(\gamma_0,..,\gamma_N)\in P_{\alpha_0}..\sigma^{-N+1}P_{\alpha_{N-1}}\}$$
otherwise;
$$G_N^\beta(\alpha_{0},..,\alpha_{N-1})=1$$
if one of the $\alpha_j$'s is of the form $i\infty$, and
$$G_N^\beta(\alpha_{0},..,\alpha_{N-1})=\frac{\beta}{2\pi}
\int_{\R^2}{\Z} e^{-\beta h_N(\gamma_0,\gamma_N)}d\gamma_0d\gamma_N$$
otherwise.

Assumption (A3) ensures us that $G_N^\beta$ is bounded, independently of $\beta$, by $B(N)$
growing subexponentially with $N$.

\begin{lem} If the constant $M$, involved in the construction of the partition $P$, is chosen large enough,
then there exists $C(\beta)\geq 0$ and, for all $N\in \N^*$, a real $\beta(N)>0$, such that : for all $k$,
and for all $\alpha_0,...,\alpha_{kN-1}$,
\begin{multline*}\left(\frac{\beta}{2\pi}
\right)^{\frac{kN+1}{2}}\int_{P_{\alpha_0}..\sigma^{-kN+1} P_{\alpha_{kN-1}}}\!\!\!\!
\psi_\beta^*(\gamma_0)e^{-\beta \sum_{i=0}^{kN-1}L(\gamma_i,\gamma_{i+1})}\psi_\beta(\gamma_{kN})
d\gamma_0...d\gamma_{kN}\\ \leq C(\beta)
\prod_{j=0}^{k-1}F_{N}(\alpha_{jN},..,\alpha_{(j+1)N-1})\prod_{l=0}^{(k-1)/2} G_N
(\alpha_{2lN},..,\alpha_{(2l+1)N-1})
(1+\petito_{\beta\rightarrow
\infty})^k
\end{multline*}
for all $\beta>\beta(N)$, and with a uniform $\petito_{\beta\rightarrow\infty}$.
\end{lem}

\begin{proof} We first note that there exists $C(\beta)>0$
such that $C(\beta)^{-1/2}\leq\psi_\beta\leq C(\beta)^{1/2}$, and $C(\beta)^{-1/2}
\leq\psi_\beta^*\leq C(\beta)^{1/2}$, because they are continuous positive $Z^d$-periodic functions.

We proceed by induction on $k$; we restrict our attention to odd $k$'s : that is, the
induction goes from $k-2$ to $k$ (the argument for even $k$
is similar -- but anyway, the aim is to let $k\To +\infty$).

Remember that we have assumed $L\geq 0$, and $L=0$ on the Aubry-Mather set -- this can be achieved by
replacing $L(\gamma_0,\gamma_1)$ by $L(\gamma_0,\gamma_1)-u(\gamma_1)+u(\gamma_0)+c(\omega)$,
$u\in S_-$.

Applying Fubini's theorem, we first estimate
the integral with respect to $\gamma_{(k-1)N+1},...,\gamma_{kN}$, when $\gamma_0,..., \gamma_{(k-1)N}$ are fixed.

If
one of the $P_{\alpha_j}$'s ($j=(k-1)N,...,kN-1$) is of the form $P_{i\infty}$, we
use Lemma 2.1.6 (b), and we get
\begin{multline*}\left(\frac{\beta}{2\pi}
\right)^{\frac{N}{2}}\int_{(\gamma_{(k-1)N},..,\gamma_{kN})\in P_{\alpha_{(k-1)N}}..\sigma^{-N+1}P_{\alpha_{kN-1}}}\!\!\!
e^{-\beta\sum_{i=(k-1)N}^{kN-1}L(\gamma_i,\gamma_{i+1})}d\gamma_{(k-1)N+1}..d\gamma_{kN}\\
\leq  1= F_N(\alpha_{(k-1)N},..,\alpha_{kN-1})G_N(\alpha_{(k-1)N},..,\alpha_{kN-1})
\end{multline*}

Otherwise, we use Lemma 2.1.6 (a), and write
\begin{multline}\left(\frac{\beta}{2\pi}
\right)^{\frac{N}{2}}\int_{(\gamma_{(k-1)N},..,\gamma_{kN})\in P_{\alpha_{(k-1)N}}..\sigma^{-N+1}P_{\alpha_{kN-1}}}
e^{-\beta\sum_{i=(k-1)N}^{kN-1}L(\gamma_i,\gamma_{i+1})}d\gamma_{(k-1)N+1}..d\gamma_{kN}
\\ \leq F_N(\alpha_{(k-1)N},..,\alpha_{kN-1})(1+\petito)
\left(\frac{\beta}{2\pi}
\right)^{\frac{1}{2}} \int_{(\gamma_{(k-1)N},\gamma_{kN})\sqin P_{\alpha_{(k-1)N}}..\sigma^{-N+1}P_{\alpha_{kN-1}}} e^{-\beta h_N(\gamma_{(k-1)N},\gamma_{kN})}
d\gamma_{kN}
\end{multline}

We now integrate
with respect to $\gamma_{(k-2)N+1},...,\gamma_{(k-1)N-1}$, while $\gamma_0,..., \gamma_{(k-2)N}$ and $\gamma_{(k-1)N}$
are still fixed.

If
one of the $\alpha_j$'s ($j=(k-2)N,...,(k-1)N-1$) is of the form $i\infty$, we
use Lemma 2.1.6 (c), and we get
\begin{multline*}\!\!\!\!\!\!\!\!\!\!\left(\!\frac{\beta}{2\pi}\!
\right)^{\frac{N-1}{2}}\!\!\!
\int_{(\gamma_{(k-2)N},..,\gamma_{(k-1)N})\in P_{\alpha_{(k-2)N}}...P_{\alpha_{(k-1)N-1}}}\!\!\!\!\!\!\!\!\!\!\!\!\!\!\!\!\!\!\!\!\!\!\!\!\!\!\!
e^{-\beta\!\sum_{i=(k\!-\!2)\!N}^{\!\!(k\!-\!1)N\!-\!1}\!\! L(\gamma_i,\!\gamma_{i+1}\!)}\!d\gamma_{\!(k\!-\!2)N\!+\!1}
..d\gamma_{\!(k\!-\!1\!)\!-\!1}\\
\leq 1\leq F_N(\alpha_{(k-2)N},..,\alpha_{(k-1)N-1})
\end{multline*}

Otherwise, we use Lemma 2.1.6 (a), and we get 

\begin{multline}\left(\!\frac{\beta}{2\pi}\!
\right)^{\frac{N-1}{2}}\!\!\!
\int_{(\gamma_{(k-2)N},..,\gamma_{(k-1)N})
\in P_{\alpha_{(k-2)N}}..\sigma^{-N+1}P_{\alpha_{(k-1)N-1}}}
e^{-\beta\sum_{i=(k-2)N}^{(k-1)N-1}L(\gamma_i,\gamma_{i+1})}d\gamma_{(k-2)N+1}..d\gamma_{(k-1)N-1}\\
\leq F_N(\alpha_{(k-2)N},..,\alpha_{(k-1)N-1})
\end{multline}
if $\beta$ is large enough. This last bound 
does not depend on $\gamma_0,..., \gamma_{(k-2)N}$

Finally, integrating with respect to $\gamma_{(k-1)N}$, and combining
the estimates (2.1.3) and (2.1.4), we have proved :
\begin{multline*}\left(\frac{\beta}{2\pi}
\right)^{\frac{kN+1}{2}}\int_{P_{\alpha_0}..\sigma^{-kN+1} P_{\alpha_{kN-1}}}\!\!\!\!
e^{-\beta \sum_{i=0}^{kN-1}L(\gamma_i,\gamma_{i+1})}
d\gamma_0...d\gamma_{kN}\\ \leq \left(\frac{\beta}{2\pi}
\right)^{\frac{(k-2)N+1}{2}}\int_{P_{\alpha_0}..\sigma^{-(k-2)N+1} P_{\alpha_{(k-2)N-1}}}\!\!\!\!
e^{-\beta \sum_{i=0}^{(k-2)N-1}L(\gamma_i,\gamma_{i+1})}
d\gamma_0...d\gamma_{(k-2)N}
F_{N}(\alpha_{(k-2)N},..,\alpha_{(k-1)N-1})
\\ \times F_{N}(\alpha_{(k-1)N},..,\alpha_{kN-1})G_N
(\alpha_{(k-1)N},..,\alpha_{kN-1})
(1+\petito_{\beta\rightarrow
\infty})
\end{multline*}

which allows to prove Lemma 2.1.6 by induction.

\end{proof}

Let us turn to the left hand side of (2.1.2), which we will try to bound below before letting $n=kN$ tend to $\infty$.
Since $\mu$ is a minimizing measure, we note that the terms are non zero only if all the $P_{\alpha_i}$'s are
included in $\{\mid \gamma_1-\gamma_0\mid\leq M\}$ (if $M$ is large enough); besides, the cylinder 
$P_{\alpha_0}.... \sigma^{-n+1} P_{\alpha_{n-1}}$ must contain a trajectory in the Mather set. By Assumptions (A1) and (A2),
if $\eps$ has been chosen small enough, and
if $(\gamma_0,\gamma_1,...,\gamma_n)$ is a configuration belonging to such a cylinder, there is exactly one
minimizer $(\gamma_0,\bar\gamma_1,...,\bar\gamma_{n-1},\gamma_n)$,
in the cylinder, of the function : $$L(\gamma_0,\xi_1)+\sum_{i=1}^{n-2}L(\xi_i,\xi_{i+1})+L(\xi_{n-1}, \gamma_n)$$ The action of this
minimizer is, by definition,
$$h_n(\gamma_0,\gamma_n)=L(\gamma_0,\bar\gamma_1)+\sum_{i=1}^{n-2}L(\bar\gamma_i,\bar\gamma_{i+1})+
L(\bar\gamma_{n-1}, \gamma_n)$$ 

As previously, we want to use the Laplace method to estimate the left hand side of (2.1.2). But since we need
to do it uniformly in the length $n$ of the path, we shall be more careful than previously.

Applying a Taylor formula to the function $L(\gamma_0,\xi_1)+\sum_{i=1}^{n-2}L(\xi_i,\xi_{i+1})+L(\xi_{n-1},
\gamma_n)$ at the minimizer $(\bar \gamma_1,\bar \gamma_2,\cdots,
\bar \gamma_{n-1})$, we can write
\begin{multline*}\int_{P_{\alpha_0}..\sigma^{-n+1} P_{\alpha_{n-1}}}
\psi^*_\beta(\gamma_0)e^{-\beta \sum_{i=0}^{n-1}L(\gamma_i,\gamma_{i+1})} \psi_\beta(\gamma_{n})d\gamma_0...d\gamma_{n}\\
=\int_{P_{\alpha_0}..\sigma^{-n+1} P_{\alpha_{n-1}}}
\psi^*_\beta(\gamma_0)e^{-\beta h_n(\gamma)-\frac{\beta}{2}\;_{N-1}\! A^{\prime\prime}(\gamma).
(\gamma-\bar \gamma)^2
-\beta R_n(\gamma)}\psi_\beta(\gamma_{n})d\gamma_0..d\gamma_{n}
\end{multline*}
where the remainder $R_n$ is given by the integral formula :
$$R_n(\gamma)=\int_0^1\frac{(1-t)^2}{2} A^{(3)}(\bar \gamma +t(\gamma-\bar \gamma)).((\gamma-\bar \gamma))^3dt$$
so that
\begin{multline*}\mid R_n(\gamma)\mid \leq C\parallel \gamma-\bar \gamma\parallel_3^3 \leq C
\parallel \gamma-\bar \gamma\parallel_\infty \parallel \gamma-\bar \gamma\parallel_2^2\leq C\eps
\parallel \gamma-\bar \gamma\parallel_2^2
\end{multline*}
since the third derivative of $L$ is bounded.

We also know (\cite{Fa2}) that $h_n$ is a Lipschitz function (with lipschitz constant independent on $n$),
and that $h_n(\gamma_0,\gamma_{n})=0$ if $(\gamma_0,\gamma_{n})$ are the endpoints of a trajectory in
the Mather set : as a consequence,
$\mid h_n(\gamma)\mid\leq C\eps$ uniformly on the cylinder
$P_{\alpha_0}..\sigma^{-n+1} P_{\alpha_{n-1}}$, as soon as it contains a trajectory in the Mather set. 

Thus, 

\begin{multline*}
\int_{P_{\alpha_0}..\sigma^{-n+1} P_{\alpha_{n-1}}}\psi_\beta^*(\gamma_0)
e^{-\beta \sum_{i=0}^{n-1}L(\gamma_i,\gamma_{i+1})}\psi_\beta(\gamma_{n})d\gamma_0...d\gamma_{n}\\
\geq C(\beta)^{-1}e^{-\beta C\eps}
\int_{P_{\alpha_0}..\sigma^{-n+1} P_{\alpha_{n-1}}}
e^{-\beta (\frac{1}{2}\;_{N-1}\! A^{\prime\prime}(\gamma)+C\eps I_{n-1}).(\gamma_1-\bar \gamma_1,..,\gamma_{n-1}-\bar \gamma_{n-1})^2}d\gamma_0...d\gamma_{n}\\
\end{multline*}

\begin{lem} For all $\alpha>0$,
there exists $r(\alpha)>0$ such that :
if $A$ is an invertible
symmetric tridiagonal matrix with $\mid A_{i, i+1}\mid \leq 1$,
then $$\parallel A^{-1}\parallel_2\leq \alpha$$
implies
$$\parallel A^{-1}\parallel_\infty \leq r(\alpha)$$
independently of the dimension.
 \end{lem}

\begin{proof} For $1\leq j\leq n$, let $f^j=A^{-1}e^j$, where $(e^j)$ is the canonical
base of $\R^n$. Note that
$$\parallel A^{-1}\parallel_\infty =\sup_k \sum_j \mid f^j_{k}\mid=\sup_j \sum_k\mid f^j_{k}\mid$$
since $A^{-1}$ is symmetric.

Let us fix $j$, and denote $f=f^j$. For $m>j$, we define a vector $^m\! f$ with coordinates
$$^m\! f_k=0$$ for $k<m$, and $$^m\! f_k=f_k$$ for $k\geq m$. Then $\eta= A.^m\! f$
has coordinates
$$\eta_{m-1}=A_{m-1,m}f_m$$
$$\eta_{m}=-A_{m,m-1}f_{m-1}$$
and
$$\eta_k=0$$otherwise.

Since, by assumption,
$$\parallel\; ^m\! f\parallel_2 \leq \alpha \parallel\eta\parallel_2,$$
we get, for all $m>j$,
$$P_m:=\sum_{k\geq m}\mid f_k\mid^2\leq
\alpha^2(\mid f_m\mid^2+\mid f_{m-1}\mid^2)
$$ 

Then the lemma p 128 of \cite{AMB} yields
$$f_k\leq \alpha^{-1}\left(\frac{2\alpha^{2}}{1+(1+4\alpha^{4})^{1/2}}\right)^{k-j/2}$$
for $k\geq j$, so that
$$\sum_{k\geq j}\mid f^j_k\mid \leq
\sum_{k=0}^{+\infty}\alpha^{-1}\left(\frac{2\alpha^{2}}{1+(1+4\alpha^{4})^{1/2}}\right)^{k/2}=:r(\alpha)/2$$
We can use a similar trick for $k<j$, and get that
$$\sum_{1\leq k\leq n}\mid f^j_k\mid\leq r(\alpha),$$ independently of $j$ and of the dimension $n$.
\end{proof}

\begin{coro} There exists $\rho(\eps)$ such that, for all $n$, for all
$\gamma\in W$,
$$\parallel \left(\;_{n-1}\! A^{\prime\prime}(\gamma)+2C\eps I_{n-1}\right)^{-1/2}\parallel_\infty \leq \frac{1}{\rho(\eps)}$$ 
\end{coro}

\begin{proof} Obviously, the spectrum of
$\;_{n-1}\! A^{\prime\prime}(\gamma)+2C\eps I_{n-1}$ is included in an interval $[\eps,\lambda]$ independent
of the dimension $n$. Let $\cal C$ be a closed contour in $\C\setminus\R^-$, going once around $[\eps,\lambda]$.
The matrix 
$\left(\;_{n-1}\! A^{\prime\prime}(\gamma)+2C\eps I_{n-1}\right)^{-1/2}$ is given by
$$\left(\;_{n-1}\! A^{\prime\prime}(\gamma)+2C\eps I_{n-1}\right)^{-1/2}=\frac{1}{2i\pi}\int_{\cal C}
z^{-1/2}\left(zI_{n-1}-(\;_{n-1}\! A^{\prime\prime}(\gamma)+2C\eps I_{n-1})\right)^{-1}dz$$
Now, for all $z\in{\cal C}$,
$$\parallel (zI_{n-1}-(\;_{n-1}\! A^{\prime\prime}(\gamma)+2C\eps I_{n-1})^{-1}\parallel_2$$
is bounded, independently of $n$, by $$\alpha(z)=\sup_{x\in[\eps,\lambda]}\frac{1}{\mid z-x\mid}$$
By Lemma ?, $$\parallel (zI_{n-1}-(\;_{n-1}\! A^{\prime\prime}(\gamma)+2C\eps I_{n-1})^{-1}\parallel_\infty\leq
r(\alpha(z))$$
independently of $n$, and
$$\parallel \left(\;_{n-1}\! A^{\prime\prime}(\gamma)+2C\eps I_{n-1}\right)^{-1/2}\parallel_\infty 
\leq\frac{1}{2\pi}\int_\gamma
\mid z\mid^{-1/2}r(\alpha(z))dz:=\frac{1}{\rho(\eps)}$$
\end{proof}

Coming back to (2.1.2), we can write
\begin{multline}\left(\frac{\beta}{2\pi}\right)^{\frac{n+1}{2}}
\int_{P_{\alpha_0}...\sigma^{-n+1} P_{\alpha_{n-1}}}\psi_\beta^*(\gamma_0)
\exp\left(-\beta \sum_{i=0}^{n-1}L(\gamma_i,\gamma_{i+1}) \right)
\psi_\beta(\gamma_{n})d\gamma_0...d\gamma_{n}\\
\geq C(\beta)^{-1}e^{-\beta C\eps}\left(\frac{\beta}{2\pi}\right)^{\frac{n+1}{2}}
\int_{P_{\alpha_0}...\sigma^{-n+1} P_{\alpha_{n-1}}}
e^{-\beta (\frac{1}{2}\;_{n-1}\! A^{\prime\prime}(\gamma)+C\eps I_{n-1}).(\gamma_1-\bar
\gamma_1,\cdots,\gamma_{n-1}-\bar \gamma_{n-1})^2
}d\gamma_0...d\gamma_{n}\\
\geq C(\beta)^{-1}e^{-\beta C\eps}\left(\frac{\beta}{2\pi}\right)^{\frac{n+1}{2}}
\int_{(\gamma_0,\gamma_n)\sqin(\cal{M}_\omega\cap
P_{\alpha_0}..\sigma^{-n+1} P_{\alpha_{n-1}})^{\eps/4}}\!\!\!\!\!\!\!\!\!\!
d\gamma_0 d\gamma_n\times \\
\int_{\parallel(\gamma_j-\bar\gamma_j)\parallel_\infty\leq\eps/4}
\exp\left(-\beta (\frac{1}{2} \;_{n-1}A^{\prime\prime}(\gamma_0,\gamma_{n})+
C\eps I_{n-1}).(\gamma_1-\bar
\gamma_1,\cdots,\gamma_{n-1}-\bar \gamma_{n-1})^2
\right)d\gamma_1..d\gamma_{n-1}\\
\geq C(\beta)^{-1}e^{-\beta C\eps}\left(\frac{\beta}{2\pi}\right)^{\frac{n+1}{2}} \int_{(\gamma_0,\gamma_n)\sqin(\cal{M}_\omega\cap
P_{\alpha_0}..\sigma^{-n+1} P_{\alpha_{n-1}})^{\eps/4}}\!\!\!\!\!\!\!\!\!
\!\!\!
d\gamma_0 d\gamma_n\times
\\
\int_{\parallel (\frac{\;_{n-1}\! A^{\prime\prime}(\gamma_0,\gamma_n)}{2}+C\eps I_n)^{1/2}.
(\gamma-\bar\gamma)\parallel_\infty\leq \rho(\eps)\eps/4}
e^{-\beta (\frac{1}{2}\;_{n-1}\!A^{\prime\prime}(\gamma_0,\gamma_{n})+C\eps I_{n-1}).(\gamma_1-\bar
\gamma_1,\cdots,\gamma_{n-1}-\bar \gamma_{n-1})^2
}d\gamma_1..d\gamma_{n-1}\\
= C(\beta)^{-1}e^{-\beta C\eps}\left(\frac{\beta}{2\pi}\right)
\int_{(\gamma_0,\gamma_n)\sqin(\cal{M}_\omega\cap
P_{\alpha_0}..\sigma^{-n+1} P_{\alpha_{n-1}})^{\eps/4}}\!\!\!\!\!\!\!\!\!\!
d\gamma_0 d\gamma_n\times
\\ \frac{1}{[\;_{n-1}\! A^{\prime\prime}(\gamma_0,\gamma_n)
+2C\eps I_{n-1}]^{1/2}}\times
\frac{1}{(2\pi)^{(n-1)/2}}\int_{\parallel y\parallel_\infty\leq \sqrt{\beta}\rho(\eps)\eps/4}
e^{-\frac{(y,y)}{2}}dy_1..dy_{n-1}\\
\geq
\frac{\eps^2}{16}\left(\frac{\beta}{2\pi}\right)
\frac{1}{\max_\alpha[\;_{n-1}\! A^{\prime\prime}(\gamma_0,\gamma_n)+2C\eps I_{n-1}]^{1/2}}
(1-e^{-\beta \rho(\eps)^2\eps^2/32})^{n-1}
\end{multline}
The $max_\alpha$ in the last line is, of course, taken over all the $(\gamma_0,\gamma_n)\sqin(\cal{M}_\omega\cap
P_{\alpha_0}..\sigma^{-n+1} P_{\alpha_{n-1}})^{\eps/4}$. 

To get the last inequality, we have used the following estimate on tails of the Gaussian distribution on $\R$ :
\begin{equation*}\frac{1}{(2\pi)^{1/2}}\int_{\mid y\mid\geq Y}
e^{-\frac{\mid y\mid^2}{2}}dy \leq 2\frac{e^{-\frac{\mid Y\mid^2}{2}}}{Y}
\end{equation*}
which yields, in dimension $n-1$,
\begin{equation}\frac{1}{(2\pi)^{(n-1)/2}}\int_{\parallel y\parallel_\infty \leq Y}
e^{-\frac{( y,y)^2}{2}}dy_1..dy_{n-1} \geq (1-e^{-\frac{( Y,Y)^2}{2}})^{n-1}
\end{equation}
for $Y>2$.

Let us summarize in a lemma what we have just proved.

\begin{lem} Let $\alpha\in W^P$ intersect the Mather set. Then, for all $(\gamma_0,\gamma_n)\sqin(\cal{M}_\omega\cap
P_{\alpha_0}..\sigma^{-n+1} P_{\alpha_{n-1}})^{\eps/4}$,
\begin{multline*}\left(\frac{\beta}{2\pi}\right)^{\frac{n-1}{2}}
\int_{P_{\alpha_0}..\sigma^{-n+1} P_{\alpha_{n-1}}}\!\!\!\!\!\!\!
\exp\left(-\beta \sum_{i=0}^{n-1}L(\gamma_i,\gamma_{i+1}) \right)
d\gamma_1...d\gamma_{n-1}\\
\geq
\frac{1}{[\;_{n-1}\! A^{\prime\prime}(\gamma_0,\gamma_n)+2C\eps I_{n-1}]^{1/2}}
(1-e^{-\beta \rho})^{n-1}
\end{multline*}
for some $\rho=\rho(\eps)^2\eps^2/32 >0$ depending only on $\eps$.
\end{lem}

To resume the proof of Theorem 2.1.1, taking $n=kN$, and putting together
Lemmas 2.1.7 and 2.1.10,
we can deduce from inequality (2.1.2) the following :
\begin{multline}-\sum \mu(P_{\alpha_0}.... \sigma^{-kN+1} P_{\alpha_{kN-1}})\log \mu(P_{\alpha_0}...
\sigma^{-kN+1} P_{\alpha_{kN-1}})
-\log( C(\beta)
\frac{\eps^2}{16}e^{-\beta C\eps})\\
+(kN-1)\log(1-e^{-\beta \rho})-\frac{1}{2}\sum_\alpha \mu(P_{\alpha_0}.. \sigma^{-kN+1} P_{\alpha_{kN-1}})
\log(\max_\alpha  [\;_{kN-1}\!A^{\prime\prime}(\gamma)+2C\eps I_{kN-1}])
\\ \leq
-\sum \mu_\beta(P_{\alpha_0}.... \sigma^{-kN+1} P_{\alpha_{kN-1}})
\log \mu_\beta(P_{\alpha_0}....\sigma^{-kN+1} P_{\alpha_{kN-1}})+
\log C(\beta)+k\log(1+\petito_{\beta\rightarrow 
\infty})\\
+k\sum \mu_\beta(P_{\alpha_0}.... \sigma^{-N+1} P_{\alpha_{N-1}})
\log F_N(\alpha_0,..,\alpha_{N-1})
+\frac{(k+1)\log B(N)}{2}
\end{multline}
for $\beta$ large enough.

We notice that $\log\max_\alpha[\;_{n-1}\!A^{\prime\prime}(\gamma)+2C\eps I_{n-1}]$, as a
function of the sequence $(\alpha_0,..,\alpha_{n-1})$,
has the following subadditivity property : if $(\alpha_0,..,\alpha_{n-1})$ intersects the Mather set,
then
\begin{multline*}\log \max_\alpha[\;_{n-1}\!A^{\prime\prime}(\gamma)+2C\eps I_{n-1}]\leq
\log \max_\alpha[\;_{m}\!A^{\prime\prime}(\gamma)+2C\eps I_{m}]
+ \log\max_\alpha[\;_{n-1-m}\!A^{\prime\prime}(\sigma^m\gamma)+2C\eps I_{n-1-m}]
\end{multline*}
This follows straightforwardly from Lemma 2.1.3.

As a consequence, if $\mu$ is an (invariant) minimizing measure, then
$$\frac{1}{kN}\sum \mu(P_{\alpha_0}.. \sigma^{-kN+1} P_{\alpha_{kN-1}})
\log\max_\alpha  [\;_{kN-1}\!A^{\prime\prime}(\gamma)+2C\eps I_{kN-1}]$$
converges to its infimum, as $k\To +\infty$. And in particular, the limit is less than
$$\frac{1}{N}\sum \mu(P_{\alpha_0}.. \sigma^{-N+1} P_{\alpha_{N-1}})
\log\max_\alpha  [\;_{N-1}\!A^{\prime\prime}(\gamma)+2C\eps I_{N-1}]$$

Thus, if we divide both side of (2.1.7) by $kN$ and let $k$ tend to $\infty$ ($\beta$ being kept
fixed), we get the inequality :
\begin{multline*}h_\sigma(\mu,P)-\frac{1}{2}\int_{W^P} \frac{1}{N}\log\max_\alpha [\;_{N-1}\!A^{\prime\prime}
(\gamma)+2C\eps
I_{N-1}]
d\mu^P(\alpha)-\log(1-e^{-\beta \rho})\\
\leq h_\sigma(\mu_\beta,P)+\int_{W^P} \frac{1}{N}\log F_N(\alpha)d\mu_\beta^P(\alpha)
+\frac{1}{N}(\petito_{\beta\rightarrow\infty}) +\frac{\log B(N)}{2N}
\end{multline*}

Now, let $\beta\To +\infty$ -- or at least, take a sequence $\beta_k$ such that $\mu_{\beta_k}$ converges
weakly to $\mu_\infty$. Supposing
that $\mu_\infty$ does not charge the boundary of the elements of the partition -- otherwise we could
always modify slightly the partition so that this assumption is satisfied -- we get
\begin{multline*}h_\sigma(\mu,P)-\frac{1}{2}\int_{W^P} \frac{1}{N}\log
\max_\alpha[\;_{N-1}\!A^{\prime\prime\max}(\gamma)
+2C\eps
I_{N-1}])
d\mu^P(\alpha)\\ \leq h_\sigma(\mu_\infty,P)+\int_{W^P} \frac{1}{N}\log F_N(\alpha)d\mu^P_\infty(\alpha)
+\frac{\log B(N)}{2N}
\end{multline*}
The point in fixing $N$ was to integrate only functions depending on a finite number of coordinates, so as to be
able to pass to the weak limit.

Now, letting $\eps\To 0$, and recalling the definition of $F_N$,
\begin{multline*}h_\sigma(\mu)-\frac{1}{2}\int \frac{1}{N}\log  [\;_{N-1}\!A^{\prime\prime}(\gamma)]
d\mu(\gamma)\\ \leq h_\sigma(\mu_\infty)-\frac{1}{2}\int \frac{1}{N}\log [\;_{N-1}\!A^{\prime\prime}(\gamma)]
d\mu_\infty(\gamma)+\frac{\log B(N)}{2N}
\end{multline*}
and, finally, letting $N\rightarrow +\infty$ (and using Assumption (A3)), we get the result.

\end{proof}
\end{proof}

This ends the proof for the discretized system.

In (2.1.5) and (2.1.6) we have proved the following fact, which will be useful in
the treatment of continuous time :
\begin{lem} For all $M\geq 0$, and for all $\eps>0$, there exists
$\rho=\rho(\eps,M)>0$ such that, for all $n$, for all $nd\times nd$ block-tridiagonal positive symmetric
matrix $A$ satisfying

-- $\mid A_{i,i+1}\mid\leq M$ for all $i$.

-- $A\geq\eps I_n$,

then
\begin{multline*}\left(\frac{\beta}{2\pi}\right)^{n/2}\int_{\parallel x\parallel_\infty\leq \eps}e^{-\beta \frac{(Ax,x)}{2}}
dx_1..dx_n\geq (1-e^{-\beta\rho})^n \left(\frac{\beta}{2\pi}\right)^{n/2}\int_{\R^n}e^{-\beta \frac{(Ax,x)}{2}}
dx_1..dx_n
\\=\frac{(1-e^{-\beta\rho})^n}{[A]^{1/2}} 
\end{multline*}
for all $\beta>0$.
\end{lem}

\subsection{Elements of the proof in continuous time}

In continuous time, the proof goes along the same lines, with
a higher degree of technicality. We will not write down the proof in its
entirety, but explain how the ideas used in discrete time can be made to work
in continuous time.

Again we treat the case $d=1$.

The proof starts as previously with the construction of a partition $\tilde P$ of $C^0([0,1], \R)$ :
$$\tilde P_{ij}=\{\gamma\in C^0([0,1], \R), \gamma_0\in [i\eps, (i+1)\eps), \gamma_1\in [j\eps, (j+1)\eps)\}$$ 
for $\mid j-i\mid <\frac{M}{\eps}$, and
$$\tilde P_{i\infty}=\{\gamma, \gamma_0\in [i\eps, (i+1)\eps), \exists j,\mid j-i\mid\geq \frac{M}{\eps},
\gamma_1\in [j\eps, (j+1)\eps)\}$$
If $\eps$ is the inverse of an integer, the partition goes to the quotient $W_{[0,1]}=C^0([0,1], \R)/\Z$,
and then gives a finite partition $P$ of $W$.

We can then write the convexity inequality (2.1.1) and try to follow the same steps.
\vspace{.3cm}

{\bf Definition of the hessian of the action, and of its determinant.}
For one moment, let us denote $\bar H$ the affine Hilbert space $H_{[0,t]}^x$ (respectively $H_{[0,t]}^{x,y}$), and $H$ its tangent
space $H_{[0,t]}^0$ (respectively $H_{[0,t]}^{0,0}$); similarly, we denote $\bar W$ the affine Banach space $W_{[0,t]}^x$ ($W_{[0,t]}^{x,y}$),
and $W$ its tangent space. Then $\bar H$ is densely immersed into $\bar W$, and $H$ is densely immersed
into $W$.

The action $\A : \bar H\To \R$ is twice differentiable, and its second derivative at a point $\gamma$,
$d^2\A(\gamma)$, is a symmetric bilinear form on $H$; one may write it as
$$d^2\A(\gamma).\xi.\xi=\langle \A^{\prime\prime}(\gamma)\xi,\xi\rangle$$
where $\A^{\prime\prime}(\gamma)$ is an autoadjoint operator on $H$ : the hessian of $\A$ at $\gamma$.

Remembering the expression of $\A$, one has $$\A^{\prime\prime}(\gamma)=I+f,$$
$f$ being defined by
$$\langle f y,y\rangle=\int_0^t V^{\prime\prime}(\gamma_s).y_s.y_s ds$$
This last bilinear form may be extended to a continuous symmetric bilinear form on $W$; and this implies
that $f$ is a trace operator (\cite{Kuo}, p.83) : the sum of the eigenvalues of $f$, $(\lambda_i)_{i\in\N}$,
is absolutely convergent.

Thus, we may define the determinant of $I+f$ as $\prod_{i\in\N}(1+\lambda_i)$, which is well defined (possibly
zero). This determinant will be non zero if and only if $-1$ is not an eigenvalue of $f$, if and only if the
operator $\A^{\prime\prime}(\gamma)$ is invertible in $H$.

If $\gamma$ is a critical point of $\A : \bar H\To\R$ such that $\A^{\prime\prime}(\gamma)$ is invertible, we will say that
$\gamma$ is a non-degenerate critical point of $\A : \bar H\To\R$.

As in the discrete time case, if $\gamma\in H_{[0,t^{\prime}]}$ for some $t^\prime\geq t$,
we
will denote $[\;_t A^{\prime\prime}(\gamma)]$
the determinant of the hessian of $A(\gamma_{|[0,t]}) : H_{[0,t]}^{\gamma_0,\gamma_t}\To\R$, at $\gamma$.

\vspace{.2cm}
{\bf Laplace method (fixed time interval).}
The analogue of Lemma 2.1.6 to continuous time can be obtained using superlinear
growth of the Lagrangian (for parts (b) and (c)); and the Laplace method for path integrals (for part (a)) :

\begin{theo} (\cite{BA}, \cite{BDS})

(a) (for Brownian motion) Assume that the action $\A : H^x_{[0,t]}\To\R$ has
only one minimum $\bar\gamma\in H_{[0,t]}^x$, which is non degenerate, and
let $V$ be a neighbourhood of $\bar\gamma$ in the uniform topology. Then
$$\int_{W^x_{[0,t]}\cap V}e^{\beta (\int_0^t V(\gamma_s)ds +\langle \omega,\gamma_t-
\gamma_0\rangle)}d{\cal W}_{[0,t]}^{\beta,x}(\gamma)
{\mathop{\sim}\limits}_{\beta+\rightarrow\infty} \frac{e^{-\beta \A(\bar\gamma)}}{[\A^{\prime\prime}(\bar\gamma)]^{1/2}}$$
where the hessian is that of $\A: H^x_{[0,t]}\To \R$ at $\bar\gamma$.

(b) (for Brownian bridge) Assume that the action $\A : H^{x,y}_{[0,t]}\To\R$ has
only one minimum $\bar\gamma$, which is non degenerate, and let $V$ be a neighbourhood of $\bar\gamma$ in the uniform topology. Then
$$\int_{W^{x,y}_{[0,t]}\cap V}e^{\beta (\int_0^t V(\gamma_s)ds +\langle \omega,\gamma_t-
\gamma_0\rangle)}d{\cal W}_{[0,t]}^{\beta,x,y}(\gamma)
{\mathop{\sim}\limits}_{\beta+\rightarrow\infty}  \frac{e^{-\beta\A(\bar\gamma)}}{
[\;_t\A^{\prime\prime}(\bar\gamma)]^{1/2}}$$
\end{theo}

These estimates are obtained, exactly as in the case of a finite
dimensional system, by applying a Taylor expansion of order $2$ of the function :
$$\gamma\mapsto \int_0^t V(\gamma_s)ds +\langle \omega,\gamma_t-
\gamma_0\rangle$$ at the minimizer $\bar\gamma$, and in the space $W^{x,y}_{[0,t]}$ (in case (b)):
\begin{multline*}\int_{W^{x,y}_{[0,t]}\cap V}e^{\beta \int_0^t V(\gamma_s)ds}d{\cal W}_{[0,t]}^{\beta,x,y}(\gamma)=\\
\int_{W^{x,y}_{[0,t]}\cap V}e^{\beta\int_0^t V^{\prime}_{\bar\gamma_s}.(\gamma_s-
\bar\gamma_s)ds +\frac{\beta}{2}\int_0^t V^{\prime\prime}_{\bar\gamma_s}.(\gamma_s-
\bar\gamma_s)^2ds +\beta R(\gamma-\bar\gamma)}d{\cal W}_{[0,t]}^{\beta,x,y}(\gamma)=
\\\int_{W^{x,y}_{[0,t]}\cap V}e^{\beta\langle\bar\gamma_s,\gamma_s-
\bar\gamma_s\rangle +\frac{\beta}{2}\int_0^t V^{\prime\prime}_{\bar\gamma_s}.(\gamma_s-
\bar\gamma_s)^2ds+\beta R(\gamma-\bar\gamma)}d{\cal W}_{[0,t]}^{\beta,x,y}(\gamma)\\
=e^{-\beta\frac{\parallel\bar\gamma\parallel^2}{2}}
\int_{W^{x,y}_{[0,t]}\cap V-\bar\gamma}e^{\frac{\beta}{2}\int_0^t V^{\prime\prime}_{\bar\gamma_s}.\gamma_s^2ds+\beta R(\gamma)}d{\cal W}_{[0,t]}^{\beta,x,y}(\gamma)
\end{multline*}
where the last line is obtained by the Cameron-Martin formula (\cite{Kuo},
p.111).

If $t=n$ and $V=P_{\alpha_0}..\sigma^{-n+1}P_{\alpha_{n-1}}$ contains
a minimizer of the action, the remainder $R(\gamma)$, given by Taylor's
integral formula, is bounded by $C\parallel\gamma\parallel_3^3$; and by
$C\eps\parallel\gamma\parallel_2^2$
on a set of relative measure $\geq (1-e^{-\beta\eps\rho})^n$.

The other ingredient is the formula
$$\int_{W^{x,y}_{[0,t]}}e^{-\beta\langle f \gamma,\gamma \rangle}d{\cal W}_{[0,t]}^{\beta,x,y}(\gamma)
=[I+f]^{-1/2}$$
valid as soon as $\langle f.,.\rangle$ is a continuous symmetric bilinear form on $H_{[0,t]}^{x,y}$
which admits a continuous extension to $W_{[0,t]}^{x,y}$.

{\bf Laplace method (lower bound, independent of the time interval).}
In order to generalize the lower bound (2.1.5) to continuous time, we are led
to check that, for all $(\gamma_0,\gamma_n)\sqin \M_\omega^\eps$,
\begin{multline}\int_{P_{\alpha_0}..\sigma^{-n+1}P_{\alpha_{n-1}}-\bar\gamma}e^{\frac{\beta}{2}
\int_0^n V^{\prime\prime}_{\bar \gamma_s}.\gamma_s^2-\beta C\eps\int_0^t
\mid \gamma_s\mid^2 ds}d{\cal W}_{[0,n]}^{\beta,\gamma_0,\gamma_n}\\
\geq (1+\petito_{\beta\rightarrow\infty})^n\int_{W^{\gamma_0,\gamma_n}_{[0,n]}} 
 e^{\frac{\beta}{2} \int_0^nV^{\prime\prime}_{\bar \gamma_s}.\gamma_s^2-\beta C\eps\int_0^t
\mid \gamma_s\mid^2 ds}d{\cal W}_{[0,n]}^{\beta,\gamma_0,\gamma_n}
\\ = (1+\petito_{\beta\rightarrow\infty})^n\frac{1}{[\;_t\A^{\prime\prime}(\bar\gamma)+2C\eps b_t]^{1/2}}
\end{multline}
where $b_t$ is the bilinear form
$\int_0^t
\mid \gamma_s\mid^2 ds$.

Given $\bar \gamma$, minimizer of $\A : H_{[0,n]}^{\gamma_0,\gamma_n}\To\R$,
let us introduce the action 
$$\tilde \A(\gamma_{|[0,t]})=\int_0^t\frac{\mid\dot\gamma_s\mid}{2}-\frac{1}{2}\int_0^t V^{\prime\prime}_{\bar \gamma_s}.\gamma_s^2+C\eps\int_0^t
\mid \gamma_s\mid^2 ds$$
for $t\leq n$, and $\gamma\in H^{\gamma_0,\gamma_n}_{[0,n]}$. Let us also introduce the function
$$Q_j(x,y)=\inf_{\gamma_j=x, \gamma_{j+1}=y}\tilde \A(\gamma_{|[j,j+1]})
$$
for $0\leq j\leq n-1$.
It is a quadratic form, as is readily seen by checking the identity of the parallelogram.

If we condition the first term of (2.2.1) with respect to $\gamma_1,...,\gamma_{n-1}$, apply
the Laplace estimates for fixed $\gamma_1,...,\gamma_{n-1}$, and then integrate
with respect to $\gamma_1,...,\gamma_{n-1}$, 
we get
\begin{multline*}\int_{P_{\alpha_0}..\sigma^{-n+1}P_{\alpha_{n-1}}-\bar\gamma}e^{\frac{\beta}{2}
\int_0^n V^{\prime\prime}_{\bar \gamma_s}.\gamma_s^2-\beta C\eps\int_0^t
\mid \gamma_s\mid^2 ds}d{\cal W}_{[0,n]}^{\beta,\gamma_0,\gamma_n}
\\
\geq (1+\petito_{\beta\rightarrow\infty})^n \int_{P_{\alpha_0}..\sigma^{-n+1}P_{\alpha_{n-1}}-\bar\gamma}
\frac{e^{-\frac{\beta}{2}(Q_0(0,\gamma_1)+...+Q_{n-1}(\gamma_{n-1},0))}}
{\prod_{j=0}^{n-1} [\tilde \A_j^{\prime\prime}]^{1/2}}d\gamma_1...d\gamma_{n-1}
\end{multline*}
where $[\tilde \A_j^{\prime\prime}]$ is the determinant of the hessian of $\tilde \A : H_{[j,j+1]}^{\gamma_j,
\gamma_{j+1}}\To\R$ at a minimum,
and does not depend on the endpoints $\gamma_j,\gamma_{j+1}$, since the action $\tilde \A$ is a
quadratic form in the path.

But now,
$$Q_0(0,\gamma_1)+Q_2(\gamma_1,\gamma_2)...+Q_{n-1}
(\gamma_{n-1},0)$$ is a quadratic form in $(\gamma_1,
...,\gamma_{n-1})\in \R^{n-1}$, which satisfies all the assumptions of Lemma
2.1.11. Also, by Assumption (A2), $P_{\alpha_0}..\sigma^{-n+1}P_{\alpha_{n-1}}-\bar\gamma$ contains a neighbourhood of $(0,..,0)$ of size $\eps$.

So, for a suitable choice of $\rho$, we can write
\begin{multline*}\int_{P_{\alpha_0}..\sigma^{-n+1}P_{\alpha_{n-1}}-\bar\gamma}
\frac{e^{-\beta(Q_0(0,\gamma_1)+...+Q_{n-1}
(\gamma_{n-1}-\bar\gamma_{n-1},0)}}
{\prod \tilde\A_j^{\prime\prime}}d\gamma_1...d\gamma_{n-1}\\
\geq (1-e^{-\beta\rho})^n \int_{\R^{n-1}}\frac{e^{-\frac{\beta}{2}(Q_0(0,\gamma_1-\bar\gamma_1)+...+Q_{n-1}
(\gamma_{n-1}-\bar\gamma_{n-1},0)}}{
\prod \tilde\A_j^{\prime\prime}}d\gamma_1...d\gamma_{n-1}\\
\geq
(1-e^{-\beta\rho})^n\int e^{-\frac{\beta}{2}
\int_0^n V^{\prime\prime}_{\bar \gamma_s}.\gamma_s^2-\beta C\eps\int_0^n
\mid \gamma_s\mid^2 ds}d{\cal W}_{[0,n]}^{\beta,\gamma_0,\gamma_n}
=\frac{(1-e^{-\beta\rho})^n}{[_n\A^{\prime\prime}(\gamma_0,\gamma_n)+2C\eps b_n]^{1/2}}\end{multline*}

Once this step has been checked, the proof goes as smoothly as in the case of discrete
time, and so far we can state the following :
\begin{prop}Let $\mu$ be an action-minimizing measure, and $\mu_\infty$ 
a limit point of $(\mu_\beta)_{\beta\To +\infty}$. Then, under assumption (A1)(A2) and (A3),
$$h_\phi(\mu)-\frac{1}{2}\int \lim_n \frac{1}{n}\log[_n\A^{\prime\prime}(\gamma)]d\mu(\gamma)
\leq h_\phi(\mu_\infty)-\frac{1}{2}\int 
\lim_n \frac{1}{n}\log[_n\A^{\prime\prime}(\gamma)]d\mu_\infty(\gamma)$$
\end{prop}

The last step is the identification of determinants of the hessian of
$\A$ with Lyapunov exponents.

{\bf Identification of determinants.}

We now prove Theorem 1.1.3.

We shall use the result obtained for discrete time systems (Lemma 2.1.5), and let the discretization step
tend to $0$, to prove the result for continuous time systems. It is sufficient to
consider the case of the time interval $[0,1]$,
from which the general case $[0,t]$ can be deduced by
a change of variables.

We recall that the Euler-Lagrange flow associated to the Lagrangian
$L(x,v)=\frac{\parallel v\parallel^2}{2}-V(x)$ is the flow on the tangent bundle $\T^d\times \R^d$,
associated to the second order equation
\begin{equation}\ddot\gamma +V^{\prime}(\gamma)=0\end{equation}
on the torus.

The equation of small variations along an orbit $\gamma$ is
\begin{equation}\ddot{y_s} +V^{\prime\prime}(\gamma_s).y_s=0
\end{equation}

To begin with, let us examine the case when the determinant of $\;_1\!\A^{\prime\prime}(\gamma)$
vanishes. This is equivalent to $\;_1\!\A^{\prime\prime}(\gamma)$
being non injective on $H^{0,0}_{[0,1]}$, and means precisely that there exists $y\in H^{0,0}_{[0,1]}$
satisfying the differential equation (2.2.3), and not vanishing identically :
$y_0=0, \dot y_0\not=0, y_1=0$. Thus, the linear map $(y_0=0,\dot y_0)\longmapsto y_1$ is not
injective, and its determinant vanishes.

Let us now consider the case when $\;_1\!\A^{\prime\prime}(\gamma)$ is invertible.

Let us divide the interval $[0,1]$ into $N$ subintervals of equal length, and
consider the following one-step discretization scheme for the equations (2.2.2) and (2.2.3) :
\begin{equation}(\Gamma_{i+1}-\Gamma_i)-(\Gamma_i-\Gamma_{i-1})+\frac{1}{N^2}V^{\prime}(\Gamma_i)=0\end{equation}
\begin{equation}(Y_{i+1}-Y_i)-(Y_i-Y_{i-1})+\frac{1}{N^2}V^{\prime\prime}(\Gamma_i).Y_i=0\end{equation}
($i=1,\cdots,N-1$).
It is nothing else than the equations of, respectively, orbits and small variations along an orbit, for
the twist diffeomorphism corresponding to the action
$$A^{N}(\Gamma_0,\Gamma_1,\cdots,\Gamma_{N-1},\Gamma_N)=N\sum_{i=0}^{N-1}
\frac{\mid \Gamma_{i+1}-\Gamma_i\mid^2}{2}-\frac{1}{N}V(\Gamma_i)$$
(starting from now, we stick to capital letters for the discretized system).

\begin{lem} There exists a constant $C$ such that, if $(\Gamma_0,
\Gamma_1,\cdots, \Gamma_N), (Y_0,\cdots,Y_N)$ are solutions
of (2.2.4) and (2.2.5), and if $\gamma(t), y(t)$ ($t\in[0,1]$) are solutions of (2.2.2),
(2.2.3) with initial conditions satisfying
\begin{eqnarray*}\mid \gamma_0-\Gamma_0\mid\leq A/N\\
\mid\gamma^{\prime}_0-\Gamma_1\mid\leq A/N\\
\mid y_0-Y_0\mid\leq A/N\\
\mid \dot y_0-Y_1\leq A/N
\end{eqnarray*} then
\begin{multline*}\max_{k=0,\cdots,N}\{\mid \gamma_{k/N}-\Gamma_k\mid,\mid \dot\gamma_{k/N}-N(
\Gamma_{k+1}-\Gamma_k)\mid,\mid y_{k/N}-Y_k\mid,\mid \dot y_{k/N}-N(Y_{k+1}-Y_k)\mid \}
\\ \leq \frac{C+A}{N}\end{multline*}
uniformly in $N$.
\end{lem}

\begin{proof} This is a straightforward application of Theorems 16.2.2
 and 16.2.3 of \cite{Scha}, applied to $$u(t)=(\gamma(t),\dot\gamma(t),
y(t),\dot y(t))$$ and to the sequence $$U_k=(\Gamma_{k-1}, N(\Gamma_k-
\Gamma_{k-1}), Y_{k-1}, N(Y_k-
Y_{k-1}))$$ obtained by the discretization scheme.
\end{proof}

The second
derivative of $A^N$ with respect to variations of $\Gamma_1,...,\Gamma_{N-1}$ takes the form
$$d^2 \;_{N-1}\! A^N(\Gamma).Y.Y=N\sum \parallel Y_{i+1}-Y_i\parallel^2 + \frac{1}{N}\sum_{i=1}^{N-1}
V^{\prime
\prime}(\Gamma_i).Y_i.Y_i$$
($Y_0=0$, $Y_N=0$).

In Lemma 2.1.5, we have precisely shown that, for any $N$, the determinant of $(Y_0=0,Y_1)\longmapsto Y_N$ is equal to the determinant
of the bilinear form $\frac{1}{N}d^2 (\,_{N-1}\! A^N)(\Gamma)$ 
with respect to the euclidean structure $\sum_i \parallel Y_i\parallel^2$. An elementary calculation shows that it is equal to
$N^d$ times the determinant of 
$\frac{1}{N}d^2(\,_{N-1}\! A^N)(\Gamma)$ with respect to the euclidean structure
$\sum_{i=0}^{N-1} \parallel Y_{i+1}-Y_i\parallel^2$. It is, equivalently,
$N^d$ times the determinant of $d^2 (\,_{N-1}\! A^N)(\Gamma)$ with respect to the scalar product 
$\langle Y,Y \rangle =N\sum_{i=0}^{N-1} \parallel Y_{i+1}-Y_i\parallel^2$.
We now stick to this euclidean structure, and consider the corresponding hessian $\;_{N-1}\! A^{N\prime\prime}(\Gamma)$.

We notice that $H_N^{0,0}:=\{(Y_0,Y_1,\cdots,Y_{N-1},Y_N), Y_0=Y_N=0\}\simeq\R^{(N-1)d}$,
endowed with the euclidean structure
$N\sum_{i=0}^{N-1} (Y_{i+1}-Y_i)^2$, can be imbedded
in the Hilbert space $H_{[0,1]}^{0,0}$ as the $(N-1)d$-dimensional subspace 
of fields $y$ which vary affinely on each $[k/N, (k+1)/N]$; an element of $H_N^{0,0}$,
seen as an element of $H_{[0,1]}^{0,0}$, is defined by the values $y_{k/N}=Y_k$.
We note that the orthogonal projection $p_N$ from $H_{[0,1]}^{0,0}$ to $H_N^{0,0}$ is precisely given by
$$y\mapsto (y_{k/N})_{k=1,\cdots,N-1}$$

In terms of operators, we can write the hessians
$$\;_1\!\A^{\prime\prime}(\gamma)=I+f $$
and $$\;_{N-1}\! A^{N\prime\prime}(\Gamma)=I+F$$
where $f$ and $F$ are defined by
$$\langle f y,y\rangle=\int_0^1V^{\prime\prime}(\gamma_t).y_t.y_t dt$$
and
$$\langle F Y,Y\rangle=\frac{1}{N}\sum_{i=1}^{N-1}V^{\prime
\prime}(\Gamma_i).Y_i.Y_i$$

We extend $\;_{N-1}\! A^{N\prime\prime}(\Gamma)$ to a symmetric operator on $H_{[0,1]}^{0,0}$, by setting $F=
F\circ p_N$. Of course, the operator $F$ depends on
$N$, but we shall neglect to show it in the notations.

We want to use the convergence of the discretization scheme (Lemma 2.2.3) to prove,
by taking the limit $N\To +\infty$, 
that the determinant of the operator $\;_1\!\A^{\prime\prime}(\gamma)$ coincides with that of the linear map $(y_0=0,
y^{\prime}_0)\longmapsto y(1)$.

It follows from Lemma 2.2.3 that, given the initial conditions $\Gamma_0=\gamma_0$, $\Gamma_1=\dot \gamma_0
$, $Y_0=y_0=0$, the determinant of the linear map $
Y_1\longmapsto \frac{Y_N}{N}$ will converge to that of $y^{\prime}_0\longmapsto
y_1$, if we let $N\To +\infty$. As we know, the former one is equal to the determinant of the
hessian $\;_{N-1}\! A^{N\prime\prime}(\Gamma)$.

So, we want to prove that the determinant of the operator
$\;_{N-1}\! A^{N\prime\prime}(\Gamma)$ converges (as $N\To\infty$) to that of $\;_1\!\A^{\prime\prime}(\gamma)$,
defined as the infinite product of eigenvalues of
$\;_1\!\A^{\prime\prime}(\gamma)$, or equivalently, $\exp(\mbox{tr } \log \;_1\!\A^{\prime\prime}(\gamma))$. 
We choose a holomorphic logarithm defined outside a half-line which does not intersect the spectrum
of $\;_1\!\A^{\prime\prime}(\gamma)$, and $\log \;_1\!\A^{\prime\prime}(\gamma)$ is defined by
\begin{equation}\frac{1}{2i\pi}\int_{\cal C}\log z . (z-\;_1\!\A^{\prime\prime}(\gamma))^{-1}dz\end{equation}
where $\cal C$ is the contour shown in Figure 1.

\begin{figure}
\begin{center}
\resizebox{10cm}{!}{\input{contour.pstex_t}}
\end{center}
\caption[e]{Definition of $\log\;_1\!\A^{\prime\prime}(\gamma)$.}
\end{figure}

We write $H_{[0,1]}^{0,0}=\overline{\cup_{N=2^n} H_N^{0,0}}$, noting that, if we take only
diadic subdivisions of the interval, the union is increasing.

On $H_N^{0,0}$ we consider the orthonormal basis consisting of functions $(e_1,\cdots,e_{N-1})$,
whose graphs are represented simultaneously (up to normalisation) on Figure 2 (in the case $d=1$, $N=16 $).

\begin{figure}
\begin{center}
\resizebox{10cm}{!}{\input{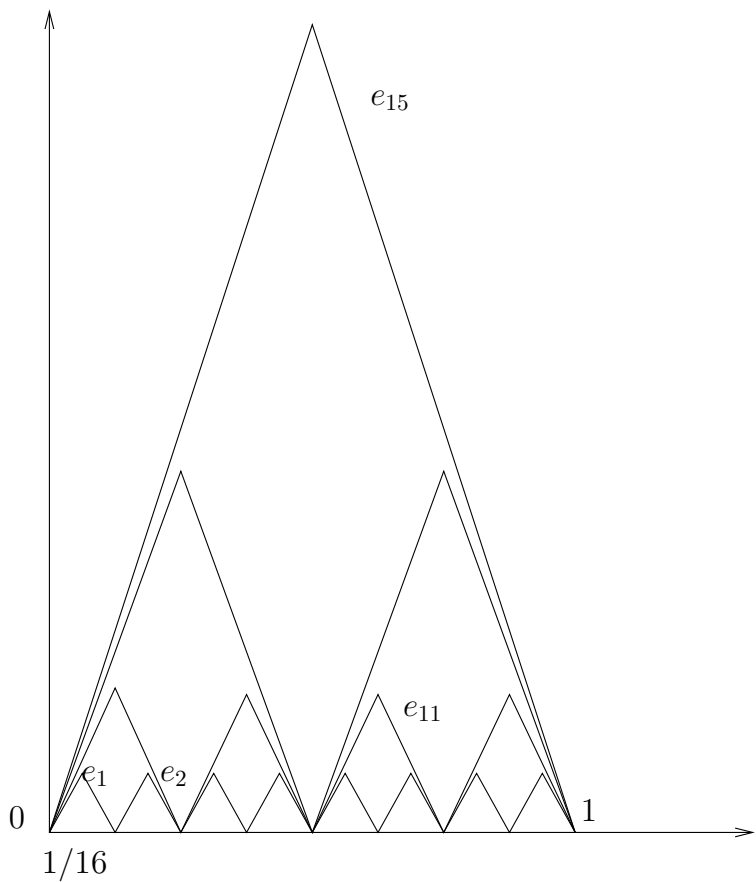}}
\end{center}
\caption[e]{An orthonormal basis of $H_N^{0,0}$.}
\end{figure}

A crucial feature is that there are $Nd/2$ elements of the basis
supported on intervals of length $2/N$, $Nd/4$ elements
supported on intervals of length $4/N$, and, more generally, $Nd/{2^l}$
elements
supported on intervals of length $2^l/N$ (for $l\leq
n=\log_2 N$. 

By definition of the trace,
\begin{eqnarray*}\mbox{tr } \log \;_1\!\A^{\prime\prime}(\gamma)=\lim_{N\To\infty}\sum_{i=1}^{N-1}\langle 
\log \;_1\!\A^{\prime\prime}(\gamma) .e_i,e_i\rangle \end{eqnarray*}

What we need to show is that this limit is the same as the limit :
$$\lim_{N\To\infty}\sum_{i=1}^{N-1}\langle 
\log \;_{N-1}\! A^{N\prime\prime}(\Gamma) e_i,e_i\rangle $$
(which we know exists).
 
Recall that $$\;_1\!\A^{\prime\prime}(\gamma)=I+f $$
and $$\;_{N-1}\! A^{N\prime\prime}(\Gamma)=I+F$$
By the expressions of
$f$ and $F$ as well as Lemma 2.2.3,
$\parallel f-F\parallel\leq \frac{C}{\sqrt N}$ (Riemann sums converge at the rate $1/\sqrt{N}$ for functions in $H_{[0,1]}$); this implies that the spectrum of
$F$ lies inside the contour $\cal C$, if $N$ is large enough. Both $\log \;_1\!\A^{\prime\prime}(\gamma)$ and
$\log \;_{N-1}\! A^{N\prime\prime}(\Gamma)$
can then be expressed thanks to a contour integral like (2.2.6). 

Thus, in order to estimate
$$\sum_{i=1}^{N-1}\langle 
\log \;_1\!\A^{\prime\prime}(\gamma) .e_i,e_i\rangle -\sum_{i=1}^{N-1}\langle 
\log \;_{N-1}\! A^{N\prime\prime}(\Gamma) e_i,e_i\rangle $$we are led to evaluate
$$\sum_{i=1}^{N-1}\langle ((z-I-f)^{-1}-(z-I-F)^{-1}).e_i,e_i\rangle$$
for all $z\in{\cal C}$.

We write
$$(z-I-f)^{-1}-(z-I-F)^{-1}=(f+I-z)^{-1}(f-F)(F+I-z)^{-1}$$
\begin{multline*}\sum_{i=1}^{N-1}\langle ((z-I-f)^{-1}-(z-I-F)^{-1}).e_i,e_i\rangle\\
=\sum_{i=1}^{N-1}\langle (f+I-z)^{-1}(f-F)(F+I-z)^{-1}e_i, e_i\rangle\\
=\sum_{i=1}^{N-1}\langle (F+I-z)^{-1}(f+I-z)^{-1}(f-F)e_i, e_i\rangle
= \sum_{i=1}^{N-1}\langle (f-F)e_i, (f+I-z)^{-1}(F+I-z)^{-1}e_i\rangle
\end{multline*}
using the property of the trace (the fact that $(F+I-z)^{-1}$ preserves $H_N^{0,0}$ is crucial),
as well as the fact that $f$ and $F$ are symmetric.

As we already mentioned, $Nd/2^l$ of the $e_i$'s vanish outside an interval
of length $2^l/N$ :

Let $\chi$ be a function in $H_N^{0,0}$ which vanishes outside an interval
$I$, and let $\zeta$ be any function in $H_{[0,1]}^{x,y}$.
Then
$$\langle f\chi,\zeta\rangle=\int_I V^{\prime\prime}(\gamma_t).\chi_t.\zeta_t dt$$
and
\begin{multline*}\langle F\chi,\zeta\rangle=
\frac{1}{N}\sum_{i/N\in I}V^{\prime\prime}(\Gamma_i).\chi_{i/N}.\zeta_{i/N}\\ =\frac{1}{N}
\sum_{i/N\in I}V^{\prime\prime}(\gamma_{i/N}
).\chi_{i/N}.\zeta_{i/N}
+\mid I\mid.\parallel \chi\parallel.\parallel \zeta\parallel O(1/N)
\end{multline*}
after Lemma 2.2.3.

Besides,
$$\int_I V^{\prime\prime}(\gamma_t).\chi_t.\zeta_t dt
-\frac{1}{N}
\sum_{i/n\in I}V^{\prime\prime}(\gamma_{i/N}
).\chi_i.\zeta_i
=\mid I\mid\parallel \chi\parallel
.\parallel \zeta\parallel O(1/\sqrt{N}),$$
a Riemann sum estimate for functions in $H_{[0,1]}$.

Applying this to the $\chi=e_i$'s and $\zeta=(f+I-z)^{-1}(F+I-z)^{-1}e_i$ and summing over $1\leq i\leq N-1$, we get
$$\mid\sum_{i=1}^{N-1}\langle ((z-I-f)^{-1}-(z-I-F)^{-1}).e_i,e_i\rangle\mid \leq C \frac{\log_2 N}{\sqrt{N}}\parallel (f+I-z)^{-1}(F+I-z)^{-1}
\parallel$$
and
$$\mid \sum_{i=1}^{N-1}\langle 
\log \;_1\!\A^{\prime\prime}(\gamma) .e_i,e_i\rangle- \langle 
\log \;_{N-1}\! A^{N\prime\prime} (\Gamma)e_i,e_i\rangle \mid \leq C \frac{\log_2 N}{\sqrt N}\int_{\cal C}\log z \parallel (f+I-z)^{-1}(F+I-z)^{-1}
\parallel dz$$
which tends to zero as $N\To +\infty$.

\end{document}